\hoffset=0.2in
\input amstex
\documentstyle{amsppt}
\magnification=1100
\nologo

\documentstyle{amsppt}

\nologo
\define\oh{{\operatorname H}}
\define\ohf#1{{\oh^{inf}_{#1}}}
\define\e{{\epsilon}}
\define\fg{{\frak g}}
\define\fa{{\frak a}}
\define\fk{{\frak k}}

\define\fd{{\frak d}}
\define\fm{{\frak m}}
\define\fn{{\frak n}}
\define\fh{{\frak h}}
\define\fc{{\frak c}}

\define\fo{{\frak o}}
\define\fp{{\frak p}}
\define\fu{{\frak u}}
\define\fq{{\frak q}}
\define\fs{{\frak s}}
\define\fl{{\frak l}}

\define\br{{\Bbb R}}
\define\bc{{\Bbb C}}
\define\AR{{A_\Bbb R}}

\define\GR{{G_\Bbb R}}
\define\HR{{H_\Bbb R}}
\define\NR{{N_\Bbb R}}
\define\SR{{S_\Bbb R}}
\define\UR{{U_\Bbb R}}
\define\KR{{K_\Bbb R}}
\define\ar{{\fa_{\Bbb R}}}
\define\pr{{\fp_{\Bbb R}}}
\define\qr{{\fq_{\Bbb R}}}
\define\gr{{\fg_{\Bbb R}}}
\define\hr{{\fh_{\Bbb R}}}
\define\kr{{\fk_{\Bbb R}}}
\define\mr{{\fm_{\Bbb R}}}
\define\nr{{\fn_{\Bbb R}}}
\define\sr{{\fs_{\Bbb R}}}

\define\ur{{\fu_{\Bbb R}}}

\define\ad{{\operatorname{ad}}}
\define\Aut{\operatorname{Aut}}

\define\BR{\Bbb R}
\define\BC{\Bbb C}

\define\Hom{\operatorname{Hom}}
\define\Mor{{\operatorname{Mor}}}
\define\Ad{\operatorname{Ad}}
\define\Ker{\operatorname{Ker}}
\define\dist{\operatorname {dist}}

\topmatter
\title On the Geometry of Nilpotent Orbits
\endtitle
\author Wilfried Schmid and Kari Vilonen
\endauthor
\address {Department of mathematics, Harvard University, Cambridge, MA
02138, USA}
\endaddress
\email schmid\@math.harvard.edu
\endemail
\address{Department of mathematics, Northwestern University, Evanston, IL 60208,
USA}
\endaddress
\email vilonen\@math.northwestern.edu
\endemail
\thanks W.Schmid was partially supported by NSF
\endgraf K.Vilonen was partially supported by NSA, NSF, the Guggenheim Foundation, and MPI Bonn
\endthanks
\endtopmatter

\document

\vskip 4\jot

\subheading{\bf 1. Introduction}
\vskip 2\jot

In this paper we describe certain geometric features of nilpotent orbits in a real semisimple Lie algebra $\gr$. Our tools are Ness' moment map \cite{N} and the proof of the Hodge-theoretic $SL_2$-orbit theorem \cite{S,CKS}; our aim is a better understanding of the Kostant-Sekiguchi correspondence \cite{Se}.

Let us recall the nature of the correspondence. We choose a Cartan decomposition $\gr = \kr \oplus \pr$, which we complexify to $\fg = \fk \oplus \fp$. Four groups will be of interest: the automorphism group $G = \Aut(\fg)^0$, the real form $\GR = \Aut (\gr)^0$, the connected subgroup $K$ with Lie algebra $\fk$, and $\KR = \GR \cap K$, which is maximal compact in both $\GR$ and $K$. Sekiguchi \cite{Se} and Kostant (unpublished) establish a bijection between the set of nilpotent $\GR$-orbits in $\gr$ on the one hand and, on the other hand, the set of nilpotent $K$-orbits in $\fp$ -- this is the Kostant-Sekiguchi correspondence.

Our proof \cite{SV2} of a representation theoretic conjecture of Barbasch and Vogan depends on a particular geometric description of the correspondence. In very rough terms, our version of the correspondence amounts to an explicit (but subtle) deformation of any nilpotent $K$-orbit in $\fp$ into the $\GR$-orbit that it corresponds to. Earlier \cite{SV2} we had reduced this result -- theorem 7.22 below -- to certain geometric statements about nilpotent orbits. These statements -- lemmas 8.5 and 8.10 -- are proved in the final section of this paper. Along the way, we obtain several results on nilpotent orbits that look interesting in their own right. What we do has implications for Kronheimer's instanton flow \cite {Kr}: the flow is real analytic at infinity, with a power series expansion that we describe recursively. 

To give some idea of our methods, we consider a nilpotent $\GR$-orbit $\Cal O$ in $\gr\nomathbreak - \nomathbreak \{0\}$. Ness' moment map \cite{N} is a real analytic, $\KR$-invariant map $m : \Bbb S(\Cal O) \to \pr$; here $\Bbb S(\Cal O) \simeq \BR^+ \backslash \Cal O$ denotes the set of unit vectors in $\Cal O$. The square norm $\|m\|^2$ assumes its minimum value exactly along a $\KR$-orbit in $\Bbb S(\Cal O)$, which we call the {\it core\/} of $\Cal O$, and denote by $C(\Cal O)$. Each point of the core determines, and is determined by, an embedding of $\fs\fl(2,\BR) \hookrightarrow \gr$, compatibly with the Cartan involutions. This fact -- in effect, a refined version of the Jacobson-Morozov theorem -- is a crucial ingredient of Sekiguchi's description of his correspondence. The core contains much information about the orbit; for example, $\Cal O$ is $\KR$-equivariantly and real analytically isomorphic to $T_{C(\Cal O)} \Cal O$, the normal bundle of the core.

The properties of nilpotent $\GR$ orbits we mentioned so far all carry over to nilpotent orbits attached to involutions: if $\HR \subset \GR$ is the fixed point group of an involutive automorphism $\sigma : \GR \to \GR$, then $\HR$ acts on the nilpotents in the $(-1)$-eigenspace of $\sigma$ on $\gr$. Orbits of this type have cores, which again can be characterized as the set of minima of $\|m\|^2$, and orbits in this setting are again isomorphic to the normal bundles along their cores. Since $K$ is the group of fixed points of the Cartan involution, this discussion applies to nilpotent $K$-orbits in $\fp$. The core of any such orbit $\Cal O_\fp$ corresponds to a $\KR$-orbit of Cartan-compatible embeddings of $\fs\fl(2,\Bbb R)$ into $\gr$, just as in the case of a nilpotent $\GR$-orbit. Orbits of the two types are Sekiguchi-related precisely when their cores coincide via the description of cores in terms of embeddings of $\fs\fl(2,\Bbb R)$ into $\gr$. This shows, in particular, that the cores of any two Sekiguchi-related orbits are $\KR$-equivariantly, real analytically isomorphic.

Not only are the cores of Sekiguchi-related orbits isomorphic, but also their normal bundles. We show this by giving a description, inspired by the nilpotent orbit theorem \cite{CKS,S}, of the fibers of the normal bundles, in terms of Cartan-compatible linear maps $\fs\fl(2,\BR) \to \gr$. Since the orbits are isomorphic to the normal bundles of the cores, we thus get $\KR$-equivariant, real analytic isomorphisms between related orbits. The existence of isomorphisms of this type had been deduced earlier from Kronheimer's results \cite{Kr} by Vergne \cite{Ve}.

The description of the normal bundles, in conjunction with arguments in \cite{CKS,S}, leads to our refinements of Kronheimer's results. We recall those results in \S3, and state and prove the refinements in \S5. Neither the logic nor the exposition of the proof of our version of the Kostant-Sekiguchi correspondence depends on these two sections. 

We wish to thank David Vogan for informative discussions. In particular, he alerted us to the fact that the isomorphism between a nilpotent orbit and the normal bundle of its core is a particular instance of a general property of homogeneous spaces of reductive Lie groups.
\vskip 3\jot

\subheading{\bf 2. Nilpotent orbits and the moment map}
\vskip 1\jot

We consider a real semisimple Lie algebra $\gr$, and let $\GR$ denote the identity component of  $\Aut (\gr)$. Further notation: $\KR\subset\GR$ is a maximal compact subgroup,
$$
\gr \ = \ \kr \oplus \pr
\tag2.1
$$
is the Cartan decomposition, and   $\theta : \gr \to \gr$ the Cartan involution. We define the inner product
$$
(\zeta_1,\zeta_2) \ = \ -B(\zeta_1,\theta\zeta_2) \qquad\qquad (\,\zeta_1,\zeta_2\in \gr\,)
\tag2.2
$$
in terms of the Killing form $B$. It is positive definite and $\KR$-invariant. We use the term ``Killing form" loosely: a $\GR$-invariant symmetric bilinear form which is negative definite on $\kr$.

Ness \cite{N} has defined a moment map for linear group actions. In our situation, it is a $\KR$-invariant, real algebraic map
$$
m\,:\, \gr-\{0\}\ \longrightarrow \ \pr \,,
\tag2.3
$$
described  implicitly by the equation
$$
(m(\zeta),\eta) \ = \ \frac 1 {2\,\|\zeta\|^2} \,\left(\frac {d\ }{dt}\,\|\Ad\exp(t\eta)\zeta\|^2\right)|_{t=0}\,.
\tag2.4
$$
As $\eta$ runs over $\gr$ in this equation, $m(\zeta)$ becomes determined as vector in $\gr$. But the inner product is $\KR$-invariant, hence $m(\zeta)$ does lie in $\pr$. The $\KR$-invariance also implies
$$
m(\Ad(k)\zeta) \ = \ \Ad(k)(m(\zeta)) \qquad (\,k\in\KR\,)\,,
\tag2.5
$$
i.e., the map $m$ is $\KR$-equivariant. To get an explicit formula for $m(\zeta)$, we calculate:
$$
\aligned
&\frac 1 {2}\left(\frac {d\ }{dt}\,\|\Ad\exp(t\eta)\zeta\|^2\right)|_{t=0} \ =  \ ([\eta,\zeta],\zeta)\ = \ -\,B([\eta,\zeta],\theta\zeta)\
\\
&\qquad= \ -\,B(\eta,[\zeta,\theta\zeta])\ = \ B(\eta,\theta{[\zeta,\theta \zeta]})\ = \ -\,(\eta,[\zeta,\theta\zeta])\ = \ -\,([\zeta,\theta\zeta],\eta)\,,
\endaligned
\tag2.6
$$
for every test vector $\eta \in \gr$, hence
$$
m(\zeta) \ = \ -\,\frac {[\zeta,\theta\zeta]}{\| \zeta\|^2}\,.
\tag2.7
$$
The moment map is invariant under scaling, hence descends to the projectivized Lie algebra $\Bbb P(\gr)$. For our purposes, it is preferable to work on
$$
\Bbb S(\gr) \ = \  \Bbb R^+\backslash(\gr-\{0\})\,,
\tag2.8a
$$
the universal (two-fold) cover of $\Bbb P(\gr) = \Bbb R^*\backslash(\gr-\{0\})$. Note that
$$
\Bbb S(\gr) \ \cong \ \{\,\zeta\in \gr \mid \|\zeta\|^2=1\,\}\,;
\tag2.8b
$$
however, to see the action of $\GR$, one must think in terms of the description (2.8a) of $\Bbb S(\gr)$.

For our next statement, we fix a particular nilpotent $\GR$-orbit $\Cal O\subset\gr-\{0\}$. By  Jacobson-Morozov,  any $\zeta\in\Cal O$ can be embedded in an essentially unique $\fs\fl_2$-triple.  In other words, there exist $\tau$, $\zeta_-$ in $\gr$ such that
$$
[\tau,\zeta]= 2\zeta\,, \ \ \ [\tau,\zeta_-]= -2\zeta_-\,, \ \ \ [\zeta,\zeta_-]= \tau\,,
\tag2.9
$$
$\tau$ is unique up to conjugacy by the centralizer  of $\zeta$ in $\GR$, and $\zeta_-$ becomes unique once $\tau$ has been chosen. In particular, the orbit $\Cal O$ determines $\tau$ up to $\GR$-conjugacy. Thus, when we re-scale $B$ by requiring
$$
B(\tau,\tau) \ = \ 2\,,
\tag2.10
$$
the normalization depends on the orbit $\Cal O$, not on the particular choice of $\zeta$. By construction, the re-scaled $B$ restricts to the linear span of $\zeta,\zeta_-,\tau$ as the trace form of  $\fs\fl(2,\Bbb R)$, to which this linear span is isomorphic. The one parameter subgroup of $\GR$ generated by $\tau$ normalizes $\zeta$ and acts on it via $\Bbb R^+$. This establishes the well-known fact that nilpotent orbits are invariant under scaling by positive numbers. The action of $\KR$ on the nilpotent orbit $\Cal O$ commutes with scaling, so the product group $\KR \times \Bbb R^+$ acts on $\Cal O$.

\proclaim{2.11 Lemma}  A point $\zeta\in\Cal O$ is a  critical point\ of the function $\zeta\mapsto \|m(\zeta)\|^2$ if and only if  there exists  a real number $a$, $a<0$, such that
$$
[[\zeta,\theta\zeta],\zeta] \ = \ a\,\zeta \ \ \ \ \text{and}\ \ \ \ [[\zeta,\theta\zeta],\theta\zeta] \ = \ -a\,\theta\zeta\,.
$$
The set of critical points is non-empty and consists of a single $\KR\times \Bbb R^+$-orbit. The function $\| m\|^2$ on $\Cal O$ assumes its minimum value exactly on the critical set.
\endproclaim

\demo{Proof} The theorem follows readily from an adaptation of \cite{N, theorems 6.1, 6.2} to the case of real group actions \cite{Ma}. It is also possible to argue directly in our particular situation, as follows. To begin with,  $\zeta$ is a critical point if and only if $\ad (m)(\zeta)$ normalizes the line $\Bbb R\zeta$; this comes down to a short calculation, as in the proof of \cite{N, theorem 6.1}.  Hence $\zeta$ is a critical point if and only if
$$
[[\zeta,\theta\zeta],\zeta] \ = \ a\,\zeta
\tag2.12a
$$
for some $a\in\Bbb R$. Applying $\theta$ to both sides, we find
$$
[[\zeta,\theta\zeta],\theta\zeta] \ = \ - a\,\theta\zeta\,.
\tag2.12b
$$
Next we argue that (2.12a), plus the nilpotency of $\zeta$, forces $a<0$. Indeed,  $[\zeta,\theta\zeta]$ lies in the (-1)-eigenspace of $\theta$, i.e., in $\pr$\,, on which $B$ is positive definite. Thus
$$
a \|\zeta\|^2 \ = \ -B([[\zeta,\theta\zeta],\zeta], \theta\zeta) \ = \ -B([\zeta,\theta\zeta], [\zeta,\theta\zeta])\ =\ -\|[\zeta,\theta\zeta]\|^2 \ \leq 0\,.
$$
Equality cannot hold: write $\zeta=\zeta_1+\zeta_2$ with $\zeta_1\in\kr$, $\zeta_2\in\pr$; $[\zeta,\theta\zeta]=0$ implies $[\zeta_1,\zeta_2]=0$; both summands are semisimple, making $\zeta$ simultaneously semisimple and nilpotent -- impossible, since $\Cal O\neq \{0\}$ by assumption. This gives the first assertion of the lemma. Continuing with the assumption that $\zeta$ is a critical point, we rescale $\zeta$ by a positive multiple to make $a=-2$. Then, if we set $\tau=-[\zeta,\theta\zeta]$ and $\zeta_-=-\theta\zeta$, the triple $\zeta,\tau,\zeta_-$ is a strictly normal S-triple in the sense of Sekiguchi \cite{Se}. The set of all $\zeta\in\Cal O$ which can be embedded into a strictly normal S-triple consists of exactly one $\KR$-orbit \cite{Se}. Thus, as claimed, the critical set in $\Cal O$ is non-empty, and $\KR\times\Bbb R^+$ acts transitively on it. The moment map (2.4) extends naturally to the complexification $\Cal O_\BC$ of $\Cal O$ -- i.e., the orbit of $\Aut(\fg)^0$ in $\fg = \BC \otimes_\BR \gr$ passing through $\Cal O$. Any critical point of $\|m\|^2 : \Cal O \to \BR_{>0}$ remains critical for the function $\|m\|^2$ on $\Cal O_\BC$. According to \cite{N,theorem 6.2}, the set of critical points of $\|m\|^2 : \Cal O_\BC \to \BR_{>0}$ coincides with the set of minima of $\|m\|^2$ on $\Cal O_\BC$. We conclude that all critical points of $\|m\|^2$ on $\Cal O$ are minima, as asserts by the lemma.
\enddemo

Let us rephrase the lemma in slightly different terms. Since $\Bbb R^+$ acts on the nilpotent orbit $\Cal O$, we can define
$$
\Bbb S(\Cal O) \ = \ \Bbb R^+\backslash\Cal O\ \cong \ \{\,\zeta\in\Cal O \mid \|\zeta\|^2 = 1 \,\}\,
\tag2.13
$$
in analogy to (2.8). We shall call
$$
C(\Cal O) \ = \ \{\,\zeta\in \Bbb S(\Cal O)\mid\text{$\zeta$ is a critical point for $\|m\|^2$}\,\}
\tag2.14
$$
the {\it core} of $\Cal O$. The core becomes a submanifold of $\Cal O$ when we identify $\Bbb S(\Cal O)$ with the set of unit vectors in $\Cal O$: in analogy to (2.8b),
$$
\text{$C(\Cal O)$ is the set of all critical points in $\Cal O$ of unit length}\,.
\tag2.15
$$
According to lemma 2.11,
$$
\aligned
\text{a)}\ \ &\text{$C(\Cal O)$ is non-empty\,,}
\\
\text{b)}\ \ &\text{$\KR$ acts transitively on $C(\Cal O)$\,, and}
\\
\text{c)}\ \ &\text{$\Bbb R^+\cdot C(\Cal O)$ is the critical set in $\Cal O$.}
\endaligned
\tag2.16
$$
The simplest example of a pair $(\gr,\kr)$ satisfying our hypotheses is $(\fs\fl(2,\Bbb R),\fs\fo(2))$. To simplify the notation, we set
$$
\gathered
\fs_\Bbb R \ = \ \fs\fl(2,\Bbb R),\ \ \text{with Cartan involution}
\\
\theta_\fs\,:\,\fs \ \to \ \fs\,,\ \ \ \  \theta_\fs(\zeta)\ = \ -\, ^t\! \zeta\,.
\endgathered
\tag2.17a
$$
The three elements
$$
e\ = \ \pmatrix 0 &1\\0 &0\endpmatrix\,,\ \ \  f\ = \ \pmatrix 0 &0\\1 &0\endpmatrix\,,\ \ \  h\ = \ \pmatrix 1 &0\\0 &-1\endpmatrix
\tag2.17b
$$
constitute a basis of $\fs_\Bbb R$ and satisfy the relations
$$
\gathered
[h,e] \ = \ 2e\,,\ \ \ [h,f] \ = \ -2f\,,\ \ \ [e,f] \ = \ h\,,
\\
\theta_\fs(e) \ = \ -f\,,\ \ \ \theta_\fs(h) \ = \ -h\,.
\endgathered
\tag2.17c
$$
Although we are interested primarily in real Lie algebras, it is useful for certain purposes to complexify. We write $\Mor (\fs,\fg)$ for the set of non-zero Lie algebra homomorphisms from $\fs = \fs\fl(2,\Bbb C)$ to the complexification $\fg = \Bbb C \otimes_\Bbb R \gr$ of $\gr$, and define
$$
\gathered
\Mor^\Bbb R(\fs,\fg) \ = \ \{\,\Phi\in\Mor (\fs,\fg)\mid \Phi \ \text{is defined over $\Bbb R$}\,\}\,,
\\
\Mor^\theta(\fs,\fg) \ = \ \{\,\Phi\in\Mor (\fs,\fg)\mid \theta\circ\Phi =\Phi\circ\theta_\fs\,\}\,,
\\
\Mor^{\Bbb R,\theta}(\fs,\fg) \ = \ \Mor^\Bbb R(\fs,\fg)\cap \Mor^\theta(\fs,\fg) \,.
\endgathered
\tag2.18
$$
Note that $\Mor^\Bbb R(\fs,\fg)$ is naturally isomorphic to $\Mor (\sr,\gr)$, the set of non-trivial morphisms between the real Lie algebras $\sr$, $\gr$. The group $\KR$ acts on $\Mor^\Bbb R(\fs,\fg)$ through the adjoint action on $\gr$\,: $(k\,\Phi)(\zeta)=_{\text{def}} \Ad k(\Phi(\zeta))$.

\proclaim{2.19 Lemma} The map $\Phi \mapsto \Phi(e)$ establishes a $\KR$-equivariant isomorphism
$$
\{\,\Phi \in \Mor^{\Bbb R,\theta} (\fs,\fg) \mid \Phi(e) \in \Cal O\,\} \ \cong \ C(\Cal O)\,.
$$
\endproclaim

\demo{Proof} Note that any $\Phi \in \Mor^{\Bbb R,\theta} (\fs,\fg)$ is uniquely determined by its value on $e$ -- cf. (2.17c). If $\zeta = \Phi(e)$ lies in the orbit $\Cal O$, it is a critical point, as follows from lemma 2.11, coupled with the relations  (2.17c); any such $\zeta$ has unit length since the normalization (2.10) of the Killing form makes $\Phi$ an isometry, relative to the trace form on $\fs_\Bbb R$. This makes the map $\Phi \mapsto \Phi(e)$ well defined and injective. It is surjective because $a = -2$ in the proof of lemma (2.11) if and only if $\|\zeta\|=1$; in that case, the triple $\zeta$, $\zeta_-$, $\tau$ defined in that proof satisfy the same relations (2.17c) as $e$, $f$, $h$. The equivariance, finally, is obvious from the definition of the action.
\enddemo

Lemma 2.19, together with 2.16, formally implies a statement that appears, in different language, in \cite{Se}: the set of nilpotent $\GR$-orbit in $\gr - \{0\}$ corresponds bijectively to the set of $\KR$-orbits in $\Mor^{\Bbb R,\theta} (\fs,\fg)$.

The inner product (2.2), normalized as in (2.10), determines a $\KR$-invariant Riemannian metric on $\Bbb S(\gr)$. We use this metric to give meaning to the gradient vector field $\nabla\|m\|^2$ on $\Bbb S(\gr)$. Note that $\gr$ acts on $\Bbb S(\gr)$ by infinitesimal translation. For $\eta\in\gr$, $\ell(\eta)$ shall denote the vector field corresponding to $\eta$. A simple calculation shows
$$
(\nabla\|m\|^2)|_\zeta \ = \ 2\ell(m(\zeta))
\tag2.20
$$
\cite{N}. In particular, the gradient vector field -- both on $\Cal O$ and on $\Bbb S(\Cal O)$ -- is tangential to $\GR$-orbits.

\proclaim{2.21 Proposition}  The function $\|m\|^2 :\Bbb S( \Cal O)\to \Bbb R$ is Bott-Morse. It assumes its minimum value on the core $C(\Cal O)$, and has no other critical points. Its gradient flow establishes a natural $\KR$-equivariant real analytic map from $\Bbb S( \Cal O)$ to the core $C(\Cal O)$ which exhibits $C(\Cal O)$ as a strong deformation retract of $\,\Bbb S( \Cal O)$.
\endproclaim

The normalization (2.10) specifies the value of $\|m\|^2$ on $C(\Cal O)$ as 2. Thus we can conclude:

\proclaim{2.22 Corollary} The family of open sets $\{\eta\in\Bbb S(\Cal O)\mid \|m\|^2(\eta)<2+\e\}$\,,\ $\e>0$\,, forms a neighborhood basis of $C(\Cal O)$\,.
\endproclaim

Since $\Bbb S(\Cal O)\!= \BR^*\backslash \Cal O$, we can combine the retraction $\Bbb S(\Cal O)\! \to \!C(\Cal O)$ with $\BR^*\! \!\to \{1\}$ to construct a retraction of $\Cal O\,$:

\proclaim{2.23 Corollary} There exists a $\KR$-equivariant, real analytic, strong deformation retraction $\Cal O \to C(\Cal O)$\,.
\endproclaim

\demo{Proof of proposition 2.21} Recall the notion of a Bott-Morse function: the critical set is a compact manifold, and the Hessian descends to a non-degenerate bilinear form on the normal bundle. Lemma 2.11 implies that $\|m\|^2$ assumes its minimum along $C(\Cal O)$ and has no critical points outside of $C(\Cal O)$, which is surely smooth and compact. Ness \cite{N, theorem 6.2} points out that the non-degeneracy is a general property of moment maps attached to linear actions of semisimple groups. This establishes all but the final  assertion. For the last assertion, let us consider the unstable set of the gradient flow of $\|m\|^2$ associated to $C(\Cal O)$, i.e., the union of the integral curves of $\nabla\|m\|^2$ emanating from $C(\Cal O)$. Because the function $\|m\|^2$ is Bott-Morse, this set is a manifold. We shall show:
$$
\gathered
\text{the unstable set of the gradient flow  associated}
\\
\text{to the critical set $C(\Cal O)$ consists of all of $\Bbb S(\Cal O)$}\,.
\endgathered
\tag2.24
$$
The existence of a retraction from $\,\Bbb S(\Cal O)\,$ to $\,C(\Cal O)\,$ will then follow. To establish (2.24), we may work on $\,\Bbb P(\Cal O)\, = \, \{\pm1\}\backslash \Bbb S(\Cal O)\,$. This allows us to complexify the situation, replacing $\,\Bbb P(\gr)\,$ by $\,\Bbb P(\fg)\,$, the projectivization of the complexified Lie algebra $\,\fg\,$, and correspondingly $\,\Bbb P(\Cal O)\,$ by $\,\Bbb P(\Cal O_\Bbb C)\,$, the appropriate orbit of $G= \Aut^0(\fg)$. The inner product on $\gr$ extends to a hermitian inner product on $\fg$, which is preserved by $\UR$, the unique maximal compact subgroup of $G$ which contains $\KR$. The definition (2.4) of the moment map carries over to the complexified setting, where it agrees with the usual (symplectic) moment map associated to the action of $\UR$ on $\Bbb P(\fg)$ \cite{N}. According to Kirwan \cite{Ki, theorem 6.18}, the stratification defined by the gradient flow on $\Bbb P(\fg)$ is $G$-invariant. In particular, it is $\GR$-invariant. Since the gradient flow on $\Bbb P(\gr)$ is tangential to the $\GR$-orbits, we can deduce that the stratification of $\Bbb P(\gr)$ defined by the gradient flow is $\GR$-invariant.
\enddemo

The assertion of corollary 2.23 can be strengthened considerably: the orbit $\Cal O$ is isomorphic to the normal bundle of its core $C(\Cal O)$. David Vogan pointed out to us that this is a particular instance of a general fact about homogeneous spaces of semisimple Lie groups. Mostow \cite{Mo, theorem 5} proved that any quotient of a semisimple group by a semisimple subgroup fibers equivariantly over an orbit of a maximal compact subgroup, with Euclidean fibers; the fibers are then necessarily the fibers of the normal bundle. The analogous statement in general case, i.e., for the quotient of a semisimple group by a closed subgroup, can be reduced to Mostow's theorem. Below we shall sketch the argument for nilpotent orbits, since we know of no statement in the literature that would imply it.

\proclaim{2.25 Proposition} There exists a $\KR$-equivariant, real analytic isomorphism $\Cal O \cong T_{C(\Cal O)} \Cal O$.
\endproclaim

\demo{Proof} We fix a point $\Phi(e)\in C(\Cal O)$ and use $\Phi\in \Mor^{\BR,\theta}(\fs,\fg)$ to identify $\fs$ with a subalgebra of $\fg$. In particular, $e\,,f\,,h$ now all lie in $\gr$, the Cartan involution maps $e$ to $-f$ and $h$ to $-h$, and $e$ lies in $C(\Cal O)$. We define
$$
\aligned
\mr \ = \ &\text{centralizer of $\,h\,$ in $\,\gr\,$},
\\
\nr \ = \ &\text{direct sum of all eigenspaces of $\,\ad\,h\,$ in $\gr$}
\\
&\text{corresponding to strictly positive eigenvalues}\,,
\\
M_\BR\ = \ &\text{centralizer of $\,h\,$ in $\,\GR\,$},
\\
\NR\ = \ &\exp \nr\,.
\endaligned
\tag2.26
$$
Then $\,\mr\oplus\nr\subset\gr\,$ is a parabolic subalgebra and $\,M_\BR\cdot \NR\,$ (semi-direct product) the corresponding parabolic subgroup of $\,\GR\,$. Since $h\in\pr$, the Cartan involution fixes $\mr$, $\KR \cap M_\BR$ is maximal compact in $M_\BR$, and
$$
\GR\ \simeq \ \KR \times_{\KR \cap M_\BR} (M_\BR \cdot \NR )\qquad\qquad\text{(fiber product)}\,.
\tag2.27
$$
The symbol $(\GR)_e$ shall denote the centralizer of $e$ in $\GR$, with the analogous convention applying also to subgroups of $\GR$ and subalgebras of $\gr$. We claim:
$$
\aligned
\text{a)}\ \ \ &(\GR)_e\ = \ (M_\BR)_e \cdot (\NR)_e\,,\ \ \ \text{and}
\\
\text{b)}\ \ \ &(\NR)_e\ = \ \exp((\nr)_e)\,.
\endaligned
\tag2.28
$$
To see this, we suppose that $\Ad g (e) = e$, and express $g$ using the decomposition (2.27) of $\GR$ and the Cartan decomposition of $M_\BR\,$:
$$
\Ad(k\,\exp \xi \,\exp \eta)\,e\ = \ e\,, \ \ \ \text{with}\ \ k\in \KR\,,\ \xi \in \mr\cap\pr\,,\ \eta\in \nr\,.
\tag2.29
$$
Then $\,\tilde e =_{\text{def}}\Ad(\exp \xi \,\exp \eta)\,e\ = \ \Ad(k^{-1})e\,$ lies in $\,C(\Cal O)\,$. Because of (2.19), the triple $\tilde e$, $\tilde f = - \theta \tilde e$, $\tilde h = [\tilde e,\tilde f]$ satisfies the same commutation relations as the triple $e$, $f$, $h$. In particular, $[[\tilde e, \tilde f], \tilde f] = -2\tilde f$. Conjugating by the inverse of $\theta (\exp \xi \,\exp \eta) = \exp(- \xi )\,\exp(\theta \eta)$, we find
$$
[\,[\,\Ad(\exp(-\theta \eta)\,\exp(2\xi)\,\exp \eta)e\,,\,f\,]\,,\,f\,] \ = \ -\,2\,f \ = \ [\,[\,e\,,\,f\,]\,,\,f\,]\,.
\tag2.30
$$
From the definition of $\nr$, one finds that $\Ad(\exp \eta) - 1_\gr$
raises $h$-weights. Similarly, $\Ad(-\theta\exp \eta) - 1_\gr$ lowers
weights, and $\Ad(\exp\xi)$ acts semisimply with strictly positive
eigenvalues, while preserving weights. We conclude: either $\exp \eta$
commutes with $e$, or else
$$
\Ad(\exp(-\theta \eta)\,\exp(2\xi)\,\exp \eta)e \ = \ \tsize\sum_\ell
\zeta_\ell
$$
is a linear combination of weight vectors $\zeta_\ell$, with $\zeta_k \neq 0$ for at least one weight $k>2$. This latter possibility is incompatible with the identity (2.30): in any finite dimensional representation of $\fs$, $f^2$ lowers weights exactly by four and is injective on all weight spaces corresponding to weights $k \geq 2$. Conclusion: $\exp \eta \in (\NR)_e$. Arguing analogously, we find that $\Ad\exp( 2\xi)e = e$, and even $\Ad(\exp \xi)e = e$ because of the nature of the action of $\Ad(\exp\xi)$. Now, in view of (2.29), $k$ must also commute with $e$. Any element of $\KR$ that commutes with $e$ must commute with $f= -\theta e$, hence with $h = [e,f]$. This puts $k \exp \xi$ into $(M_\BR)_e$, as asserted by (2.28a). Finally, if $\,\exp \eta\,$, with $\eta \in \nr$, centralizes $e$, then so do $\,\eta = \log(\exp \eta)$ and the one parameter group generated by $\eta$. This implies (2.28b).

The centralizers of $e$ in $\KR$ and $M_\BR$ commute with all of $\fs$. In the case of $\KR$, we just gave the argument; for $M_\BR$ it follows from the observation that any two members of an $\fs\fl_2$-triple -- in our case, $e$ and $h$ -- determine the third. For emphasis,
$$
(\KR)_e\ = \ (\KR)_\fs\,,\ \ \ (M_\BR)_e\ = \ (M_\BR)_\fs\,.
\tag2.31
$$
In particular, $(\KR)_e$ and $(M_\BR)_e$ normalize both $\nr$ and $(\nr)_e$. We can choose a linear complement $\fc_\BR$ to $(\nr)_e$ in $\nr$, which is $(\KR)_e$-invariant and  $(M_\BR)_e$-invariant: we decompose $\gr$ $\sr$-isotypically; in the isotypic subspace of highest weight $r$, we take the sum of all eigenspaces corresponding to eigenvalues strictly between 0 and $r$; then $\fc_\BR$, the sum of all of these spaces for $r>0$, has the required properties. Since $\nr = \fc_\BR \oplus (\nr)_e$ (direct sum of vector spaces), 
$$
\fc_\BR \times (\nr)_e \ @>{\ \sim\ }>> \ \NR\,,\ \ \ (\zeta,\eta)\ \mapsto \ \exp \zeta \cdot \exp \eta\,,
\tag2.32
$$
is a $(\KR)_e$-invariant, $(M_\BR)_e$-invariant, real analytic isomorphism. Indeed, the diffeomorphism statement can be reduced to an assertion about nilpotent matrix groups, which can be verified using Engel's theorem; the invariance properties are a consequence of the particular choice of $\fc_\BR$. Because of (2.27-28) and (2.31-32),
$$
\aligned
\GR\ &\simeq \ \KR \times_{\KR \cap M_\BR} (M_\BR \cdot \NR)
\\
&\simeq  \ \KR \times_{\KR \cap M_\BR} M_\BR \times_{(M_\BR)_\fs} ((M_\BR)_\fs \times \fc_\BR \times (\NR)_e)
\\
&\simeq  \ \KR \times_{\KR \cap M_\BR} (M_\BR \times \fc_\BR )\times_{(M_\BR)_\fs} ((M_\BR)_\fs \times (\NR)_e)
\\
&=  \ \KR \times_{\KR \cap M_\BR} (M_\BR \times \fc_\BR )\times_{(M_\BR)_\fs} (\GR)_e\,,
\endaligned
\tag2.33
$$
as real analytic manifold with left $\KR$- and right $(\GR)_e$-action; here $(M_\BR)_\fs$ acts on $M_\BR$ by right translation and on $\fc_\BR$ by conjugation. According to \cite{Mo,theorem 5}, there exists an isomorphism  
$$
M_\BR\ \simeq \ (\KR \cap M_\BR) \times_{(\KR)_\fs} (\pr \cap \mr \cap (\mr)_\fs^\perp) \times (M_\BR)_\fs 
\tag2.34
$$
of real analytic manifolds with left $(\KR\cap M_\BR)$- and right $(M_\BR)_\fs$-action. Mostow states his decomposition theorem for connected, semisimple groups; the extension to our situation is straightforward. In the decomposition (2.34), $(M_\BR)_\fs$ and $(\KR)_\fs$ act on $\pr \cap \mr \cap (\mr)_\fs^\perp$ by conjugation. We conclude:
$$
\Cal O \ \simeq \ \GR/(\GR)_e \ \simeq \ \KR \times_{(\KR)_\fs} \left((\pr \cap \mr \cap (\mr)_\fs^\perp)\cap \fc_\BR \right)\,.
\tag2.35
$$
This is equivalent to the statement of the proposition.
\enddemo

\vskip 3\jot

\subheading{\bf 3. The instanton flow}
\vskip 1\jot

In the previous section, we described a flow on a nilpotent orbit $\Cal O$ which retracts the orbit to its core. Kronheimer has constructed a different flow, which also retracts the nilpotent orbit to its core \cite{Kr}. Let us describe his construction in slightly different language.

We continue with the notation and hypotheses of \S 2. While we are interested in a nilpotent orbit $\Cal O$ of the real group $\GR$ in the real Lie algebra $\gr$, we will work also with the complexified group $G$, the complexified Lie algebra $\fg$, and the complexification $\fs=\fs\fl(2,\Bbb C)$ of $\fs_\Bbb R =\fs\fl(2,\Bbb R)$. In analogy to (2.18), we define
$$
\gathered
\Hom (\fs,\fg) \ = \ \text{vector space of $\Bbb C$-linear maps $\Phi:\fs\to\fg$}\,,
\\
\Hom^\Bbb R(\fs,\fg) \ = \ \{\,\Phi\in \Hom (\fs,\fg)\mid \Phi \ \text{is defined over $\Bbb R$}\,\}\,,
\\
\Hom^\theta(\fs,\fg) \ = \ \{\,\Phi\in \Hom (\fs,\fg)\mid \theta\circ\Phi =\Phi\circ\theta_\fs\,\}\,,
\\
\Hom^{\Bbb R,\theta}(\fs,\fg) \ = \ \Hom^\Bbb R(\fs,\fg)\cap \Hom^\theta(\fs,\fg)
\,.
\endgathered
\tag3.1
$$
The Lie bracket can be viewed as a $G$-equivariant linear map $\wedge^2\fg\to\fg$. In the case of $\fs$, this is an isomorphism for dimension reasons, hence can be inverted to an $SL(2,\Bbb C)$-equivariant linear map $\fs\to\wedge^2\fs$. Combining the two maps, we get a symmetric bilinear pairing
$$
Q\ :\ \Hom (\fs,\fg)\otimes \Hom (\fs,\fg) \ \longrightarrow \ \Hom (\fs,\fg)\,,
\tag3.2a
$$
which is uniquely characterized by the equation
$$
Q(\Phi_1,\Phi_2)[u,v] \ = \ \frac 1 2 \left(\,[\Phi_1(u),\Phi_2(v)]-[\Phi_1(v),\Phi_2(u)]\,\right) \qquad (\,u,v\in\fs\,)\,.
\tag3.2b
$$
Note that
$$
Q(\Phi,\Phi) \ = \ \Phi \ \ \Longleftrightarrow \ \ \Phi\in\Mor (\fs,\fg)\,;
\tag3.3
$$
here, as in the previous section, $\Mor(\fs,\fg)$ denotes the set of Lie algebra homomorphisms. The pairing is defined over $\Bbb R$, i.e.,
$$
Q\ :\ \Hom^\Bbb R(\fs,\fg)\otimes \Hom^\Bbb R(\fs,\fg) \ \longrightarrow \ \Hom^\Bbb R(\fs,\fg)\,,
\tag3.4
$$
and it is compatible with the Cartan involutions, in the sense that
$$
Q\ :\ \Hom^\theta(\fs,\fg)\otimes \Hom^\theta(\fs,\fg) \ \longrightarrow \ \Hom^\theta(\fs,\fg)\,.
\tag3.5
$$
These properties are immediate consequences of (3.2b).

\proclaim{3.6 Notation}  $\Cal M$ is the set of $C^\infty$-maps $\Phi: (0,\infty)\to \Hom^{\Bbb R,\theta}(\fs,\fg)$ satisfying the three conditions
\roster
\item "a)" $\frac d {dt}\Phi(t)\ =\ -\,Q(\Phi(t),\Phi(t))$\,,\vskip 1\jot
\item "b)" $\Phi$ extends continuously to $[0,\infty)$\,,\vskip 1\jot
\item "c)" $\lim_{t\to\infty} (t\Phi(t))$ exists and lies in $\Mor^{\Bbb R,\theta} (\fs,\fg)$\,.
\endroster
For $\Phi_0\in \Mor^{\Bbb R,\theta} (\fs,\fg)$, we set $\,\Cal M(\Phi_0)\, = \, \{\,\Phi\in\Cal M\mid \tsize \lim_{t\to\infty} (t\Phi(t)) = \Phi_0\,\}\,.$ If $C(\Cal O)$ is the core of a nilpotent $\GR$-orbit $\Cal O\subset\gr$, $\Cal M(C(\Cal O))$ will denote the union of the $\Cal M(\Phi_0)$ corresponding to
morphisms $\Phi_0$ whose image $\Phi_0(e)$ under the isomorphism (2.19) lies in $C(\Cal O)$.
\endproclaim

The conditions b),c) in this definition can be restated in equivalent, but seemingly weaker form -- see below.

\proclaim{3.7 Theorem} (Kronheimer, \cite{Kr}) The space $\Cal M$ has a natural structure of $C^\infty$ manifold. Via the map $\Phi(\,\cdot\,)\mapsto \Phi(0)(e)$, this manifold is $\KR$-equi\-variantly diffeomorphic to the nilpotent orbit $\Cal O$\,.
\endproclaim

Strictly speaking, Kronheimer states this result for complex groups. Vergne \cite{Ve} observed that the statement about real groups formally follows from the result about complex groups by restriction. Kronheimer deduces the manifold structure from general properties of moduli spaces for instantons. The manifold structure also becomes apparent from our results in \S5.

To make the transition to Kronheimer's formulation, we attach to each $\Phi\in\Cal M$ a triple of $\gr$-valued functions by evaluating $\Phi(t)$ on the triple (2.17b),
$$
E(t) \ = \ \Phi(t)(e)\,, \ \ \ F(t) \ = \ \Phi(t)(f)\,, \ \ \ H(t) \ = \ \Phi(t)(h)\,.
\tag3.8a
$$
This triple completely determines the function $\Phi$. The requirement that the values $\Phi(t)$ be compatible with the Cartan involution translates into the condition
$$
F(t) \ = \ -\,\theta E(t)\,, \ \ \ H(t) \ = \ -\,\theta H(t)\,.
\tag3.8b
$$
Let us transcribe the conditions a),b),c) in the definition 3.6. The differential equation (3.6a) becomes
$$
\gathered
2E'(t) \ = \ -\,[H(t),E(t)]\,, \ \ \ 2F'(t) \ = \ [H(t),F(t)]\,, \\  H'(t) \ = \ -\,[E(t),F(t)]\,;
\endgathered
\tag3.9a
$$
the first of these follows from $Q(\Phi,\Phi)(e)=\frac 1 2 Q(\Phi,\Phi)[h,e] =\frac 1 2 [\Phi(h),\Phi(e)]$, and similarly for the others. Next,
$$
E(t),F(t), H(t) \ \  \text{extend continuously to }\ [0,\infty)\,,
\tag3.9b
$$
and finally,
$$
\gathered
\text{the limits} \ \ E_0\, = \, \lim_{t\to\infty} (tE(t))\,, \ \ F_0\, = \, \lim_{t\to\infty} (tF(t))\,, \ \ H_0\, = \, \lim_{t\to\infty} (tH(t))
\\
\text{exist and satisfy} \ \ 2E_0\,=\,[H_0,E_0]\,, \ \ 2F_0\,=\,-[H_0,F_0]\,, \ \ H_0\,=\,[E_0,F_0]\,.
\endgathered
\tag3.9c
$$
In terms of the triple, the map $\Phi\mapsto \Phi(0)(e)$ reduces to evaluating $E(t)$ at zero.

Kronheimer, who works in the context of complex nilpotent orbits, uses a triple of $\fg$-valued functions corresponding to a different basis of $\fs$. Also,
he uses the coordinate $x= -\log t$ on $\Bbb R$, which gives a slightly different appearance to the differential equation (3.9a) and the ``evaluations" $E(t)\rightsquigarrow E_0$ and $E(t) \rightsquigarrow E(0)$.

The $\gr$-valued function $\frac 1 2 H(t)$ is the logarithmic derivative of a $C^\infty$ function $g(t)$ with values in $\GR$ -- in other words, $2g(t)^{-1}g'(t)=H(t)$. Since
$$
\gathered
\frac d {dt}\left(\Ad g(t)(E(t))\right)\ =\ \Ad g(t)\left([g(t)^{-1}g'(t),E(t)]+E'(t)\right)\ =
\\
\Ad g(t)\left([\tsize\frac 1 2 H(t),E(t)]+E'(t)\right)\ = \ 0\,,
\endgathered
\tag3.10
$$
the curve $E(t)$, for $0<t<\infty$, stays inside a nilpotent $\GR$-orbit $\Cal O$. The fact that $E(0)$ and $E_0=\lim_{t\to\infty} (tE(t))$ lie in the same orbit $\Cal O$ is a consequence of Kronheimer's theorem. Because of (3.8b) and (3.9c) -- equivalently, because
$\Phi_0$ belongs to $\Mor^{\Bbb R,\theta} (\fs,\fg)$ -- $E_0$ lies in the core $C(\Cal O)$.  In particular, then, $E(0)\mapsto E_0$ exhibits $C(\Cal O)$ as the strong deformation retract of $\Cal O$. Via the isomorphism (3.7),
$$
C(\Cal O) \ \cong \ \{\,\Phi_0\in \Mor^{\Bbb R,\theta} (\fs,\fg) \mid \Phi_0(e)\in\Cal O\,\}
\tag3.11a
$$
corresponds to the  $\Hom^{\Bbb R,\theta}(\fs,\fg)$-valued functions
$$
t \ \mapsto \ \Phi(t)\ = _{\text{def}}\ \Phi_0\,(1+t)^{-1}\,,
\tag3.11b
$$
which satisfies the differential equation (3.6a) and takes the value $\Phi_0$ at $t=0$. There are two simple operations on $\Cal M(\Phi_0)$ as defined in (3.6): for $a\in\Bbb R^+$
$$
\{\,t\mapsto\Phi(t)\,\} \longrightarrow \{\,t\mapsto a\Phi(at)\,\}\,,
\tag3.12
$$
which corresponds to scaling on $\Cal O$ under  the isomorphism (3.7), and
$$
\{\,t\mapsto\Phi(t)\,\} \longrightarrow \{\,t\mapsto a\Phi(a(t+1)-1)\,\}\,,
\tag3.13
$$
$1\leq a < \infty$, which induces the homotopy between the identity map $1_\Cal O$ and the retraction $\Cal O \to C(\Cal O)$; note that (3.13) does act trivially on the functions (3.11b).

The instanton flow is a flow in $\Hom^{\Bbb R,\theta}(\fs,\fg)$, the gradient flow of the function $\Phi\mapsto\|\Phi\|^2$ on $\Hom^{\Bbb R,\theta}(\fs,\fg)$ \cite{Kr}. Via the isomorphism $\Cal M(\Cal O) \cong \Cal O$, it corresponds to the retraction (3.13), which is not a gradient flow of a function on $\Cal O$ or $\Bbb S(\Cal O)$, nor even the flow of a (time independent!) vector field. Curiously, the retraction is induced by a vector field on certain submanifolds of nilpotent orbits, namely those which arise from variations of Hodge structure \cite{S}.

The functions $\Phi(\,\cdot\,)\in \Cal M$ are real analytic: for any $t_0\in (0,\infty)$, the coefficients of the Taylor series of $\Phi(t)$ at $t=t_0$ are polynomials in $\Phi(t_0)$ by repeated differentiation of the equation (3.6a), and the radius of convergence of this Taylor series can be bounded from below by a uniform multiple of $\|\Phi(t_0)\|^{-1}$. In particular, the condition (3.6b) can be replaced by the formally weaker condition
$$
\|\Phi(t)\| \ \ \text{is bounded on} \ \ (0,\infty)\,,
\tag3.14
$$
as long as the remaining conditions are maintained. It implies the stronger condition
$$
\|\Phi(t)\| \ \ \text{extends real analytically to} \ \ [0,\infty)\,.
\tag3.15
$$
In \S 5 we shall show that the $\Phi(t)$ are real analytic even at infinity, as functions of the variable $t^{-\frac 1 2 }$.

\vskip 3\jot

\subheading{\bf 4. The normal bundle of the core}
\vskip 1\jot

The core $C(\Cal O)$ of a nilpotent $\GR$-orbit $\Cal O\subset\gr$ is a $\KR$-orbit. This fact gives the normal bundle $T_{C(\Cal O)}\Cal O$ the structure of $\KR$-homogenous vector bundle. As such, it is associated to the representation of
$$
(\KR)_\zeta \ = \ \text{isotropy subgroup at $\zeta$}\,,
\tag4.1
$$
for any particular $\zeta\in C(\Cal O)$, on the quotient
$$
[\zeta,\gr]/[\zeta,\kr] \ \cong \ (T_{C(\Cal O)}\Cal O)_\zeta\,.
\tag4.2
$$
In this section, we shall construct a $(\KR)_\zeta$-invariant linear complement to $[\zeta,\kr]$ in $[\zeta,\gr]$. We shall need this construction in subsequent sections.

We identify the base point $\zeta$ with the morphism $\Phi_0\in \Mor^{\Bbb R,\theta} (\fs,\fg)$ which corresponds to $\zeta$ via the isomorphism (2.19). To simplify the discussion, we use $\Phi_0$ to identify $\fs_\Bbb R$ with a subalgebra of $\gr$. This physically puts the generators (2.17) into
$\gr$, with $\zeta=e$. For emphasis,
$$
e,f,h\in\gr\,, \ \ \ e=\zeta\,, \ \ \ \theta e=-f\,, \ \ \ \theta h=-h\,.
\tag4.3
$$
Since $e=\zeta$ and $\theta e=-f$ generate $\fs$,
$$
(\KR)_\zeta \ \ \text{centralizes}\ \ \fs\,.
\tag4.4
$$
The commutation relations of the triple $e,f,h$ imply that $h$ acts semisimply with integral eigenvalues in any finite dimensional representation of $\fs$. Irreducible finite dimensional representations of $\fs$ are uniquely characterized by their highest $h$-weight, which can be any non-negative integer; the irreducible representation of highest weight $r$ has dimension $r+1$. We set
$$
\gathered
\fg(r) \ = \ \text{$\fs$-isotypic subspace of $\fg$ of highest weight $r$}\ =
\\
\text{linear span of all $\fs$-irreducible subspaces of heighest weight $r$}\,;
\\
\fg(r,\ell) \ = \ \text{$\ell$-weight space of $h$ in $\fg(r)$}\,.
\endgathered
\tag4.5a
$$
The irreducible $\fs$-module of highest weight $r$ has $h$-weights $r,r-2,\dots, -r$, hence
$$
\fg\ \ = \ \ \tsize \bigoplus_{r\geq 0}\  \fg(r) \ \ = \ \ \bigoplus_{r\geq 0}\ \ \bigoplus \Sb -r\leq \ell\leq r \\ \ell\equiv  r  \!\mod 2\endSb \ \fg(r,\ell)\,.
\tag4.5b
$$
The first of these decompositions is $\fs$-invariant, $\theta$-stable, and defined over $\Bbb R$.

Recall the notation (3.1). Because of (4.4), $\Hom (\fs,\fg)$ contains the $\Hom (\fs,\fg(r))$ as $(\KR)_\zeta$-invariant subspaces -- invariant with respect to the trivial action on $\fs$ and the natural one on $\fg$. Note that $\fs$ has three natural actions on $\Hom (\fs,\fg)$: via the action on $\fs$, via its embedding in $\fg$, and diagonally. The decomposition
$$
\Hom (\fs,\fg)\ \ = \ \ \tsize \bigoplus_{r\geq 0}\ \Hom (\fs,\fg(r))
\tag4.6
$$
is $\fs$-invariant with respect to all three actions, $(\KR)_\zeta$-invariant, $\theta$-stable, and defined over $\Bbb R$. The summand corresponding to $r$ is $\fs$-isotypic of highest weight  2 with respect to the first $\fs$-action, and of highest weight $r$ with respect to the second action. Thus
$$
\aligned
&\Hom (\fs,\fg(r))\ = \
\\
&\ \ \ \ \ \ \ \Hom (\fs,\fg(r))(r-2)\,\tsize\bigoplus\,\Hom (\fs,\fg(r))(r) \,\bigoplus\,\Hom (\fs,\fg(r))(r+2)\,,
\endaligned
\tag4.7
$$
with the outer index referring to the $\fs$-type with respect to the diagonal $\fs$-action. This decomposition is also $(\KR)_\zeta$-invariant, $\theta$-stable, and defined over $\Bbb R$. Note that
$$
\Hom (\fs,\fg(r))(r-2) \ = \ 0 \ \ \ \text{unless}  \ \ r \geq 2 \,.
\tag4.8
$$
We write $\Hom^{\theta}(\fs,\fg(r))$ for the intersection of $\Hom (\fs,\fg(r))$ with $\Hom^{\theta}(\fs,\fg)$, and
analogously in the case of the summands in (4.7). Our next statement describes the fiber of the normal bundle $T_{C(\Cal O)}\Cal O$ at $\zeta$ as $(\KR)_\zeta$-module.

\proclaim{4.9 Proposition} The map $\Hom (\fs,\fg) \owns \Phi\mapsto \Phi(\zeta)=\Phi(e)$ is injective on
$$
\frak d(\Phi_0) \ \ =_{\text{def}} \ \ \tsize \bigoplus_{r\geq 2}\ \Hom^{\theta}(\fs,\fg(r))(r-2)\,.
$$
The image $\frak d(\Phi_0)(\zeta)$ of $\frak d(\Phi_0)$ under this map is a $(\KR)_\zeta$-invariant linear complement to $[\zeta,\fk]$ in $[\zeta,\fg]$, and is defined over $\Bbb R$. 
\endproclaim

Let $\frak d_\Bbb R(\Phi_0)=\frak d(\Phi_0)\cap \Hom^\Bbb R(\fs,\fg)$ denote the space of real points in $\frak d(\Phi_0)$. Then $[\zeta,\gr] = \frak d_\Bbb R(\Phi_0)(\zeta)\oplus[\zeta,\kr]$, and this identifies $(T_{C(\Cal O)}\Cal O)_\zeta \cong [\zeta,\gr]/[\zeta,\kr]$ with  $\frak d_\Bbb R(\Phi_0)(\zeta)$ as $(\KR)_\zeta$-module.

\demo{Proof of 4.9}  The evaluation map $\Phi\mapsto \Phi(\zeta)=\Phi(e)$ sends $\Hom^{\theta}(\fs,\fg(r))(r-2)$ to $\fg(r)$. We can therefore argue one summand at a time. The decompositions (4.6) and (4.7) are defined over $\Bbb R$, and $\zeta=e$ is real. This reduces the problem to showing
$$
\{\,\Phi(e) \mid \Phi\in \Hom^{\theta}(\fs,\fg(r))(r-2)\,\} \ \subset \ \ad(e)\fg(r)\,,
\tag4.10a
$$
i.e., the image of the evaluation map lies in the image of $\ad (\zeta)$, and
$$
\gathered
\text{for each $\xi\in\fg(r)$, there exist $\Phi\in \Hom^{\theta}(\fs,\fg(r))(r-2)$ and $\eta\in\fk$ so that}
\\
\Phi(e) \ = \ [e,\xi+\eta]\,; \ \ \text{in this situation,}\ \ [e,\xi] \ \ \text{uniquely determines} \ \ \Phi\,.
\endgathered
\tag4.10b
$$
For the first assertion, note that $\Phi$, which is $(r-2)$-isotypic relative to the diagonal action, has components only in the $h$-weight spaces corresponding to weights between $2-r$ and $r-2$. The evaluation map is $\fs$-equivariant and $e$ has weight two, so $\Phi(e)$ cannot have a non-zero component in the $(-r)$-weight space. In particular, this forces $\Phi(e)$ to lie in the image of $\ad (e)$.

We write $\xi=\xi_\fk + \xi_\fp$ with $\xi_\fk\in\fk$, $ \xi_\fp\in\fp$, and combine $\xi_\fk$ with $\eta$. This transforms (4.10b) into the equivalent assertion
$$
\gathered
\text{for each $\xi\in\fp\cap\fg(r)$, there exist $\Phi\in \Hom^{\theta}(\fs,\fg(r))(r-2)$ and $\eta\in\fk$}
\\
\text{so that}\ \Phi(e) = [e,\xi\!+\!\eta]\,; \ \text{in this situation,}\ [e,\xi] \ \text{uniquely determines} \ \Phi\,.
\endgathered
\tag4.11
$$
The Casimir operator of $\fs$,
$$
\Omega \ = \ 2ef+2fe+h^2
\tag4.12
$$
acts by the scalar $k^2+2k$ on any $k$-isotypic $\fs$-module. For $A$ in the universal enveloping algebra of $\fs$, we let $A\Phi$  denote the effect
of $A$ on $\Phi$, acting via the diagonal $\fs$-action  on $\Hom^{\theta}(\fs,\fg(r))$; $A\circ\Phi$ and $\Phi\circ A$ shall denote the composition of $\Phi$ with the action of $A$ on, respectively, the values and arguments. Then $\xi\Phi = \ad(\xi)\circ\Phi - \Phi\circ\ad(\xi)$ if $A = \xi\in\fs$, hence
$$
\aligned
&\Omega \Phi \ = \ \Omega\circ\Phi + \Phi\circ\Omega +
\\
&\ \ \ -4\ad(e)\circ\Phi\circ\ad(f) - 4\ad(f)\circ\Phi\circ\ad(e) - 2\ad(h)\circ\Phi\circ\ad(h)\,.
\endaligned
\tag4.13
$$
Since $(r-2)^2+2(r-2)-[r^2+2r]-[2^2+2\times2]= -4(r+2)$,
$$
\aligned
&\text{for}\ \ \Phi\in\Hom^{\theta}(\fs,\fg(r)) \ \ \text{the following two conditions are equivalent:}
\\
&\text{a)}\ \ \Phi\in \Hom^{\theta}(\fs,\fg(r))(r-2)\,,
\\
&\text{b)}\ (r\!+\!2)\Phi =  \ad(e)\circ\Phi\circ\ad(f)\! +\! \ad(f)\circ\Phi\circ\ad(e)\! +\!\frac 1 2\ad(h)\circ\Phi\circ\ad(h)\,.
\endaligned
\tag4.14
$$
To construct a particular $\Phi\in \Hom^{\theta}(\fs,\fg(r))(r-2)$ amounts to specifying $\Phi(e)$ and $\Phi(h)$ in $\fg(r)$, subject to the following conditions. First, $\Phi(h)$ must lie in $\fp$ since $h\in\fp$, and secondly, the identity (4.14b) must hold when evaluated on either $e$ or $h$. The $\theta$-equivariance of $\Phi$ then forces $\Phi(f)= -\theta\Phi(e)$. The validity of (4.14b) applied to $f$ is automatic since $\Omega$ commutes with $\theta$.

Let us suppose that $\Phi\in \Hom^{\theta}(\fs,\fg(r))(r-2)$, $\xi\in\fp\cap\fg(r)$, and $\eta\in\fk$ are given subject to the condition in (4.11), i.e.,
$$
\Phi(e)\ = \ \ad(e)(\xi+\eta)\,.
\tag4.15a
$$
This implies, and is implied by,
$$
\Phi(f)\ = \ \ad(f)(-\xi+\eta)\,,
\tag4.15b
$$
and furthermore, implies
$$
\ad(h)\eta \ = \ \ad(e)(\Phi(f)) - \ad(f)(\Phi(e)) - (\ad(e)\ad(f) + \ad(f)\ad(e))\xi \,.
\tag4.15c
$$
These identities allow us to express $\Omega\eta$ in terms of $\Phi(e),\Phi(f)$, $\xi$ and the action of $\fs$. Since $\Phi(e),\Phi(f)$, $\xi$ lie in $\fg(r)$ by assumption, $\Omega\eta$ also lies in $\fg(r)$. Thus $\eta=\eta_0+\eta_1$ with $\eta_0\in\fg(0) = \ker(\Omega)$ and $\eta_1\in\fg(r)$. Both $\ad(e)$ and $\ad(f)$ annihilate $\eta_0$, so we may as well suppose that $\eta=\eta_1\in\fg(r)$. For $r=0$, the right hand sides of (4.15a,b) vanish, and $\Hom^{\theta}(\fs,\fg(r))(r-2)=0$, which means that there is nothing to prove. Thus we may assume
$$
\eta\in\fk\cap\fg(r)\,, \ \ \ \text{and}\ \ r>0\,.
\tag4.16
$$
From these hypotheses, we shall conclude
$$
\aligned
&\text{a)} \ \ \Phi(h) \ = \ -r\xi + \tsize\frac 1 {r+2}\ad(h)^2\xi + \frac 2 {r+2} \ad(h)\eta\,,
\\
&\text{b)}\ \ \eta \ = \ \frac 1 r \,\ad(h)\xi\,.
\endaligned
\tag4.17
$$
That, in turn, will imply (4.11).

To establish (4.17), we  evaluate (4.14b) on $h$ and use the commutation relations of $e,f,h$, as well as the identities (4.15a,b):
$$
\gathered
(r+2)\Phi(h) \ =  2\,\ad(e)(\Phi(f)) - 2\,\ad(f)(\Phi(e)) \ =
\\
2\,\ad(e)\ad(f)(-\xi + \eta) \ -  \ 2\,\ad(f)\ad(e)(\xi + \eta)\ = 
\\
-\,2\,(\ad(e)\ad(f)+\ad(e)\ad(f))(\xi)\ + \ 2(\ad(e)\ad(f)-\ad(e)\ad(f))(\eta)\ =
\\
-\,(\Omega-\ad(h)^2)(\xi) \ + \ 2 \,\ad(h)(\eta)\ =
\\
-\,(r^2+2r)\xi \ + \ \ad(h)^2\xi \ + \ 2 \,\ad(h)(\eta)\,,
\endgathered
\tag4.18
$$
which is the identity (4.17a). Next, we evaluate (4.14b) on $e$. We use (4.15a), the commutation relations of $e,f,h$, and (4.18):
$$
\gathered
(r+2)\,\ad(e)(\xi+\eta)\ = \ (r+2)\,\Phi(e) \ =  \ - \,\ad(e)(\Phi(h)) +\,\ad(h)(\Phi(e)) \ =
\\
- \,\ad(e)(-r\xi + \tsize\frac 1 {r+2}\,\ad(h)^2\xi + \frac 2 {r+2}\,\ad(h)\eta)\, + \,\ad(h)\ad(e)(\xi+\eta) \,.
\endgathered
\tag4.19a
$$
Note that $\ad(h)\ad(e)= \ad(e)\ad(h) +2 \,\ad(e)$, hence
$$
\aligned
&(r+2)\,\ad(e)(\xi+\eta)\ =
\\
&\ad(e)\left(r\xi -\tsize\frac 1 {r+2}\,\ad(h)^2\xi - \frac 2 {r+2}\,\ad(h)\eta+2(\xi+\eta)+\ad(h)(\xi+\eta)\right) \,.
\endaligned
\tag4.19b
$$
We bring all terms to the left and multiply through with $r+2$, to conclude
$$
\left(\ad(h)^2-(r+2)\ad(h)\right)\xi \ + \ \left((r^2+2r)-r\,\ad(h)\right)\eta\ \in\ \ker(\ad(e))\,.
\tag4.20
$$
Recall the decomposition (4.5) and write $\xi_\ell$, $\eta_\ell$ for the components of $\xi,\eta$ in $\fg(r,\ell)$. Since $h\in\fp$, the Cartan involution interchanges $\fg(r,\ell)$ and $\fg(r,-\ell)$, and
$$
\eta_{-\ell} \ = \ -\theta(\eta_\ell)\,, \ \ \ \xi_{-s} \ = \ \theta(\xi_\ell)\,,
\tag4.21
$$
since $\ad(h)\eta\in\fp$ and $\ad(h)\xi\in\fk$. The kernel of $\ad(e)$ on $\fg(r)$ is $\fg(r,r)$. This makes (4.20) equivalent to
$$
\left((r^2+2r)-r\,\ell\right)\eta_\ell \ = \ -\,\left(\ell^2-(r+2)\ell\right)\xi_\ell\ \ \ \text{if}\ \ \ell\neq r\,.
\tag4.22
$$
The same identity for $\ell=r$ follows from the case $\ell=-r$ and (4.21). Also, $\ell$ lies between $r$ and $-r$, so (4.22) is equivalent to $r\,\eta_\ell = \ell \,\xi_\ell$ for all $\ell$. That is the assertion (4.17b).

On $\fg(r)\cap \fp$, $r>0$, $\ad(e)$ is injective, so $[e,\xi]$ determines $\xi$. From (4.17a,b) and the original hypothesis (4.15a), we now conclude that $[e,\xi]$ completely determines $\Phi(h)$, $\Phi(e)$, and $\Phi(f) = -\theta\Phi(e)$. Thus $\Phi$ is indeed uniquely determined. As was pointed out earlier, the existence of $\Phi$ comes down to knowing that $\Phi(h)$ and $\Phi(e)$ lie in $\fg(r)$, that $\Phi(h)\in\fp$, and that (4.14b) is satisfied when both sides are evaluated on $h$ and $e$, respectively. The expression (4.17b) specifies $\eta$ as element of $\fg(r)\cap \fk$. Since $h\in\fp$, $\ad(h)$ interchanges $\fk$ and $\fp$. Thus (4.17a) exhibits $\Phi(h)$ as lying in $\fg(r)\cap\fp$, as required. The containment $\Phi(e)\in \fg(r)$ follows from (4.15a) and (4.16). The validity of (4.14b) when evaluated on $h$ and $e$ amounts to the two identities (4.18) and (4.19); both hold by construction. This gives us the existence and uniqueness of $\Phi$ -- in other words, the validity of (4.11). 
\enddemo

\vskip 3\jot

\subheading{\bf 5. The instanton flow at infinity}
\vskip 1\jot

In this section we use the proof of the $Sl_2$-orbit theorem of Hodge theory \cite{S,CKS} to show that the flow lines of the instanton flow are real analytic at infinity. In effect, the proof of the $Sl_2$-orbit theorem produces a real analytic isomorphism between a neighborhood of the core $C(\Cal O)$ in a nilpotent orbit $\Cal O$ and a neighborhood of the zero section in the normal bundle $T_{C(\Cal O)}\Cal O$. This isomorphism is closely related to Kronheimer's flow. We shall freely use the notation of the earlier sections, in particular that of section 4.

The decompositions (4.5-7) depend on the embedding $\Phi_0 : \sr \hookrightarrow \gr$ given by any particular choice of $\Phi_0 \in \Mor^{\Bbb R,\theta} (\fs,\fg)$; the morphism $\Phi_0$, in turn, was assumed to correspond to some $\zeta \in C(\Cal O)$ via the isomorphism (2.19). We shall now let $\zeta$ vary over the core, and cor\-respondingly $\Phi_0$ over the inverse image of $C(\Cal O)$. In this way, we tacitly regard the decompositions (4.5-7) as depending on $\zeta$, without putting this dependence into the notation. Recall the definition (3.2) of the pairing $Q$. We shall need the notion of a $Q$-polynomial: a function
$$
\gathered
(\Phi_1,\Phi_2,\dots,\Phi_k)\ \mapsto \ P(\Phi_1,\Phi_2,\dots,\Phi_k)\in \Hom^\theta(\fs,\fg)\,,
\\
\text{with arguments}\ \ \Phi_1,\dots\,,\Phi_k \in \Hom^\theta(\fs,\fg)\,,
\endgathered
\tag5.1
$$
expressible as a finite $\Bbb C$-linear combination of monomials in the $\Phi_\ell$, with $Q$ serving as ``multiplication". Note that a real $Q$-polynomial -- i.e., a $Q$-polynomial with real coefficients -- takes values in $\Hom^{\Bbb R,\theta}(\fs,\fg)$ whenever the arguments lie in this real subspace. 

\proclaim{5.2 Theorem} Every function $\Phi(t)$ in $\Cal M$ has a convergent series expansion around $\infty$, in powers of $t^{-\frac 12}$. Specifically,
$$
\Phi(t) \ = \ \Phi_0\,t^{-1}\,+\, \tsize \sum_{k\geq 2} \Phi_k\,t^{-1-\frac k2} \qquad (\, t \gg 0\,)\,,
\tag "\qquad a)"
$$
with $\Phi_0 \in \Mor^{\Bbb R,\theta} (\fs,\fg)$ and $\Phi_k \in \Hom^{\Bbb R,\theta}(\fs,\fg)$ for $k\geq 2$; there exist universal \footnote{i.e., not depending on $\Phi_0$ nor even on $\gr$, provided the dimension of $\gr$ is bounded.} $Q$-polynomials with rational coefficients $P_k(\dots)$, $k \geq 2$, such that
$$
\aligned
&\Phi_k\ = \ \Phi_k^k + P_k(\Phi_0,\Phi_2^2, \Phi_3^3, \dots,\Phi_{k-2}^{k-2})\,,
\\
&\Phi_\ell^\ell \in \Hom^{\Bbb R,\theta}(\fs,\fg(\ell))(\ell-2) \qquad (\,\ell \geq 2\,) \,,
\\
&P_k(\Phi_0,\Phi_2^2, \Phi_3^3, \dots,\Phi_{k-2}^{k-2}) \in \tsize \bigoplus_{\ell \leq k-2,\, \ell \equiv k\,(2)}\, \Hom^{\Bbb R,\theta}(\fs,\fg(\ell))\,.
\endaligned
\tag "\qquad b)"
$$
The polynomial $P_k$ is weighted homogeneous of weight $k$ when the variable $\Phi_\ell^\ell$ is given weight $\ell$, and $\Phi_0$ weight 0. Conversely, any series of this form has a positive radius of convergence, and the resulting $\Hom^{\Bbb R,\theta}(\fs,\fg)$-valued function $\Phi(t)$ satisfies the differential equation (3.6a).
\endproclaim

Loosely paraphrased, the space $\bigoplus_{\ell \geq 2} \Hom^{\Bbb R,\theta}(\fs,\fg(\ell))(\ell-2)$ parameterizes all functions $\Phi(t)$ defined for large positive values of $t$ which satisfy the differential equation (3.6a) and the limiting condition (3.6c). In the preceeding section, we had identified this direct sum with the fiber of the normal bundle $T_{C(\Cal O)}\Cal O$ at $\zeta$ when the leading coefficient $\Phi_0 \in \Mor^{\Bbb R,\theta}(\fs,\fg)$ corresponds to $\zeta \in C(\Cal O)$. Note that the power $t^{-\frac 32}$ gets skipped in the expansion (5.2) -- this reflects (4.8).

We shall verify the theorem together with the following companion statement. For any collection of data $\Phi_0 \in \Mor^{\Bbb R,\theta}(\fs,\fg)$, $\Phi_\ell^\ell \in \Hom^{\Bbb R,\theta}(\fs,\fg(\ell))(\ell-2)$ for $\ell \geq 2$, and $t>0$, let
$$
\Phi(\Phi_0,\Phi_2^2,\dots,\Phi_\ell^\ell,\dots;t)
\tag5.3
$$
denote the value of the function $\Phi(t)$ in (5.2a,b), provided the analytic continuation of the series is defined at $t$.

\proclaim{5.4 Theorem} The assignment $(\Phi_0,\Phi_2^2,\dots,\Phi_\ell^\ell,\dots) \mapsto \Phi(\Phi_0,\Phi_2^2,\dots,\Phi_\ell^\ell,\dots;1)(e)$ induces a well defined, real analytic, $\KR$-equivariant isomorphism
$$
F\ : \ \Cal U \ \longrightarrow \ F(\Cal U)\,,
$$
between a connected open neighborhood $\,\Cal U$ of the zero section in the normal bundle $T_{C(\Cal O)}\Cal O$ and $F(\Cal U)$, an open neighborhood of the core $C(\Cal O)$ in the nilpotent orbit $\Cal O$.
\endproclaim

\demo{Proof of 5.2 and 5.4} We appeal to the results of \cite{CKS, \S6}, specifically (6.8-24); these results already appear in \cite{S}, in somewhat different language. To begin with,
$$
\gathered
\text{a formal $\Hom(\fs,\fg)$-valued series} \ \ \Phi(t) = \Phi_0\,t^{-1}\,+\, \tsize \sum_{k>0} \Phi_k\,t^{-1-\frac k2}
\\
\text{is a formal solution of the differential equation}\ \ \tsize\frac{d\Phi}{dt} = -Q(\Phi,\Phi)
\\
\text{if and only if the coefficients}\ \ \Phi_k \ \ \text{can be expressed as in (5.2b),}
\endgathered
\tag5.5
$$
with certain specific $Q$-polynomials $P_k(\dots)$. The terminology ``$Q$-polynomial" is not used in \cite{CKS}. Rather, the arguments there show that the differential equation translates into the conditions
$$
\Phi_k\ =\ \Phi_k^k + A \left( \tsize\sum_{0<\ell<k}\,Q(\Phi_\ell,\Phi_{k-\ell})\right)\qquad (\,k>0\,)\,,
\tag5.6
$$
with $\Phi_k^k \in \Hom(\fs,\fg(k))(k-2)$, and with $A$ denoting a particular rational linear combination of projections to the various eigenspaces of the linear map
$$
\Hom(\fs,\fg)\ \owns \ T \ \longmapsto \ Q(\Phi_0,T)\,.
\tag5.7
$$
The eigenvalues of this linear map are rational \cite{CKS, (6.14)}, so each projection can be expressed as a rational linear combination of its powers. The coefficient $\Phi_1$ vanishes because of (4.8) and (5.6). Thus, in (5.6), the range of summation is really $2\leq \ell \leq k-2$. Now, arguing inductively, one finds rational $Q$-polynomials $P_k(\dots)$, weighted homogeneous of degree $k$ when weights are assigned as in theorem 5.2, such that
$$
\Phi_k \ = \ \Phi_k^k + P_k(\Phi_0,\Phi_2^2, \Phi_3^3, \dots,\Phi_{k-2}^{k-2})\qquad (\,k>0\,)\,.
\tag5.8
$$
The linear map (5.7) preserves the subspaces $\Hom(\fs,\fg(\ell))$, and
$$
\aligned
&Q\left(\Hom(\fs,\fg(k_1)),\Hom(\fs,\fg(k_2))\right)\ \ \subset 
\\
&\qquad\qquad \subset \ \ \tsize \bigoplus_{0\leq \ell\leq k_1+k_2,\, \ell\equiv k_1+k_2 (2)}\, \Hom(\fs,\fg(\ell))
\endaligned
\tag5.9
$$
\cite{CKS, (6.21)}. Hence
$$
P_k(\Phi_0,\Phi_2^2, \Phi_3^3, \dots,\Phi_{k-2}^{k-2}) \in \tsize \bigoplus_{\ell \leq k-2,\, \ell \equiv k\,(2)}\, \Hom(\fs,\fg(\ell))\,,
\tag5.10
$$
again by induction on $k$. This completes the verification of (5.5). For future reference, we note that
$$
\gathered
\text{for fixed $\Phi_0$, as function of $\,\Phi_2^2, \Phi_3^3, \dots,\Phi_{k-2}^{k-2}\,$ alone,}
\\
\text{$P_k(\Phi_0,\Phi_2^2, \Phi_3^3, \dots,\Phi_{k-2}^{k-2})$ has no linear and no constant term},
\endgathered
\tag5.11
$$
as follows from the homogeneity property of $P_k(\dots)$.

The construction of the $P_k$ readily implies a bound on their size: with $\Phi_0$ kept fixed, there exists a positive constant $C$ such that
$$
\|P_k(\Phi_0,\Phi_2^2, \Phi_3^3, \dots,\Phi_{k-2}^{k-2})\| \ \leq \ C^k \left( \tsize \max_{k\geq 2} \|\Phi_k^k\| \right)^k
\tag5.12
$$
\cite{CKS, (6.24)}. That, in turn, implies
$$
\aligned
\text{a)}\ \ &\text{the series}\ \ \Phi(t)\ = \ \Phi_0\,t^{-1}\,+\, \tsize \sum_{k>0} \Phi_k\,t^{-1-\frac k2}
\\
&\text{converges if}\ \ t^{\frac 12} > C^{-1} \left( \tsize \max_{k\geq 2} \|\Phi_k^k\| \right)^{-1}\,;\vspace{-2\jot}
\\
\text{b)}\ \ &(\Phi_2^2,\dots,\Phi_\ell^\ell,\dots) \mapsto \Phi(\Phi_0,\Phi_2^2,\dots,\Phi_\ell^\ell,\dots;1)\ \, \text{is a well defined, analytic}\vspace{-3.5\jot}
\\
&\text{map on some neighborhood of $\ 0\ $ in}\ \ \tsize \bigoplus_{k\geq 2}\, \Hom(\fs,\fg(k))(k-2)\,.
\endaligned
\tag5.13
$$
In particular, the map in the statement of theorem 5.4 is well defined, real analytic when restricted to a small neighborhood of 0 in any fiber $(T_{C(\Cal O)}\Cal O)_\zeta$ of the normal bundle; Kronheimer's theorem 3.7 implies that the map takes values in $\Cal O$. Because of (4.9), (5.8) and (5.11), this map sends any sufficiently small neighborhood of 0 in $(T_{C(\Cal O)}\Cal O)_\zeta$ isomorphically to a real analytic submanifold of $\Cal O$, normal to $C(\Cal O)$ at $\zeta$. The definition (3.2) exhibits $Q$ as $\KR$-equivariant pairing. We conclude that $Q$-polynomials are $\KR$-equivariant as functions of their arguments. The map $F$ is therefore both $(\KR)_\zeta$-equivariant on the fiber at $\zeta$ and globally $\KR$-equivariant. Since $\KR$ is compact and acts transitively on $C(\Cal O)$, $F$ has the properties asserted by theorem (5.4): $\KR$-equivariant, real analytic, and real analytically invertible from a small neighborhood of the zero section in $T_{C(\Cal O)}\Cal O$ to some neighborhood of $C(\Cal O)$ in $\Cal O$.

We now consider a particular curve $\Psi(t)$ in $\Cal M(\Phi_0)$. Then, if $a\geq 1$, the curve $C_a\Psi$, defined by
$$
(C_a\Psi)(t) \ = \ a\Psi(a(t+1)-1)\,,
\tag5.14a
$$
satisfies the three conditions (3.6a-c), with the same limiting morphism $\Phi_0$. In other words, $C_a\Psi \in \Cal M(\Phi_0)$, hence
$$
C_a \ : \ \Cal M(\Phi_0)\ \longrightarrow \ \Cal M(\Phi_0)\qquad (\,a \geq 1\,)\,.
\tag5.14b
$$
The condition (3.6c) implies
$$
\lim_{a\to\infty} C_a\Psi(t)\ = \ \Phi_0\,,
\tag5.14c
$$
for any fixed $t \geq 0$ -- recall: $C_a\Psi$, like every curve satisfying (3.6a-c), extends real analytically to a neighborhood of 0. In particular, for $a$ sufficiently large,
$$
(C_a\Psi)(0)(e) \in \Cal U \,,
\tag5.15
$$
with $\,\Cal U$ having the same meaning as in the statement of theorem 5.4. Since this theorem has already been proved, there exist $\Phi_k^k \in \Hom^{\Bbb R,\theta}(\fs,\fg(k))(k-2)$, $k\geq 2$, so that
$$
\Phi(t)\ =_{\text{def}} \ \Phi_0\,t^{-1}\,+\, \tsize \sum_{k\geq 2}\left(\Phi_k^k + P_k(\Phi_0,\Phi_2^2,\dots, \Phi_{k-2}^{k-2})\right)t^{-1-\frac k2}
\tag5.16a
$$
converges for $t \gg 0$, extends real analytically to $[1,\infty)$, and satisfies
$$
\Phi(1)(e) \ = \ (C_a\Psi)(0)(e)\,.
\tag5.16b
$$
Because of (5.5), $t\mapsto\Phi(t+1)$ belongs to $\Cal M(\Phi_0)$. By construction, this curve has the same image under the Kronheiner isomorphism as $C_a\Psi$, and thus must coincide with $C_a\Psi$. We conclude that $\Psi$ can be obtained from $\Phi$ by a linear coordinate change, and that $\Psi(t)$ has a convergent series expansion around infinity, in powers of $t^{-\frac12}$. Appealing once more to (5.5) and subsequent statements, we conclude that $\Psi(t)$ has the properties asserted in theorem 5.2.
\enddemo

Theorem 5.4 can be strengthened, as follows. Recall (5.14). A short calculation shows that $C_a \circ C_b = C_{ab}$, hence
$$
a \ \mapsto \ C_a \ \ \text{induces an action of the multiplicative semigroup $\Bbb R_{\geq 1}$}\,,
\tag5.17
$$
both on $\Cal M(\Phi_0)$ and on $\Cal O\cong \Cal M(C(\Cal O))$. Using the identification (4.9), we define
$$
\gathered
D_a \ : \ T_{C(\Cal O)}\Cal O \ \longrightarrow \ T_{C(\Cal O)}\Cal O\,,
\\
D_a(\Phi_0)\ = \ \Phi_0\,,\qquad D_a(\Phi_k^k)\ = \ a^{-\frac k2}\Phi_k^k\,, 
\endgathered
\tag5.18a
$$
where $\Phi_k^k\in\Hom^{\Bbb R,\theta}(\fs,\fg(k))(k-2)$. This makes sense for all $a\neq 0\,$; moreover,
$$
a \ \mapsto \ D_a \ \ \text{defines an action of the multiplicative group $\Bbb R^*$}\,,
\tag5.18b
$$
as can be checked by direct calculation. The map $F$ defined in theorem 5.4 is $\Bbb R_{\geq 1}$-equivariant with respect to the two actions (5.17-18):

\proclaim{5.19 Lemma}  For all $a \geq 1$\,,\ \ $F\circ D_a\ = \ C_a \circ F$\,.
\endproclaim

\demo{Proof} If $(\Phi_0,\Phi_2^2,\dots,\Phi_\ell^\ell,\dots)$ corresponds to a point in the domain of $F$, the series (5.16a) converges for large $t$, extends real analytically to $[1,\infty)$, and
$$
\gathered
F(\Phi_0,\Phi_2^2,\dots,\Phi_\ell^\ell,\dots) \ = \ \Phi(1)(e)\ =\ \tilde\Phi(0)(e)\,,
\\
\text{where}\ \ \tilde\Phi(t)\ =_{\text{def}} \ \Phi(t+1)\,.
\endgathered
\tag5.20
$$
The curve $\tilde\Phi(t)$ then belongs to $\Cal M(\Phi_0)$ and corresponds to $F(\Phi_0,\Phi_2^2,\dots,\Phi_\ell^\ell,\dots)$ via the Kronheimer isomorphism (3.7). Hence
$$
C_a(F(\Phi_0,\Phi_2^2,\dots,\Phi_\ell^\ell,\dots))\, = \, (C_a\tilde\Phi)(0)(e)\, = \, a\,\tilde\Phi(a-1)(e)\, = \, a\,\Phi(a)(e)\,.
\tag5.21
$$
On the other hand,
$$
(F \circ D_a)(\Phi_0,\Phi_2^2,\dots,\Phi_\ell^\ell,\dots)\ = \ F(\Phi_0,a^{-1}\Phi_2^2,\dots,a^{-\frac \ell 2}\Phi_\ell^\ell,\dots)\,.
\tag5.22
$$
When the $\Phi_\ell^\ell$ in (5.16a) are replaced by $a^{-\frac \ell 2}\,\Phi_\ell^\ell$, the series $\Phi(t)$ gets transformed into $a\,\Phi(at)$ -- this follows from the homogeneity property of the $P_k(\dots)$. Thus
$$
F(\Phi_0,a^{-1}\Phi_2^2,\dots,a^{-\frac \ell 2}\Phi_\ell^\ell,\dots)\ = \ \left( a\,\Phi(at)|_{t=1}\right) (e)\ = \ a\,\Phi(a)(e)\,,
\tag5.23
$$
completing the proof of the lemma.
\enddemo

\proclaim{5.24 Corollary} The map $F^{-1}$ extends to a real analytic isomorphism between the entire nilpotent orbit $\Cal O$ and an open neighborhood of the zero section in $T_{C(\Cal O)}\Cal O$.
\endproclaim

\demo{Proof} Given $\xi \in \Cal O$, we choose $a>1$ so large that $C_a\xi$ lies in the domain of $F^{-1}$, and set
$$
F^{-1}(\xi)\ = \ D_{a^{-1}}(F^{-1}(C_a\xi))\,.
\tag5.25
$$
This extension of $F^{-1}$ is well defined and one-to-one by (5.17-19), and real analytic by construction. It is also locally invertible, again because of (5.17-19), hence an open map with real analytic inverse.
\enddemo

\vskip 3\jot

\subheading{\bf 6. Complex groups and symmetric pairs}
\vskip 1\jot

In preparation for the next section, in which we discuss the Sekiguchi corres\-pondence, we shall restate the earlier results for symmetric pairs and complex Lie algebras.

There is not so much to say about the complex case -- complex Lie algebras can be regarded as real Lie algebras, after all. Let $\fg$ be a complex semisimple Lie algebra. As a matter of general notational convention, we set
$$
\fg^\Bbb R \ = \ \fg\,,\ \ \text{taken as Lie algebra over $\Bbb R$}\,.
\tag6.1
$$
In the situation of interest to us, $\fg$ will arise as the complexification of a real semisimple Lie algebra $\gr$, in which a Cartan decomposition $\gr = \kr\oplus\pr$ has been specified. The subalgebra
$$
\ur \ = \ \kr \tsize\oplus i\,\pr
\tag6.2
$$
is then a compact real form in $\fg^\Bbb R$, and
$$
\fg^\Bbb R \ = \ \ur \tsize\oplus i\,\ur
\tag6.3
$$
the Cartan decomposition determined by $\ur$. Further notation:
$$
\tau \, : \, \fg \ \rightarrow \fg\ \ \ \text{is complex conjugation with respect to}\ \ \ur\,.
\tag6.4
$$
In view of (6.3),
$$
\text{$\tau$\ \,is the Cartan involution on\ \,$\fg^\br$}\,.
\tag6.5
$$
We normalize the Killing forms on $\fg$, $\fg^\Bbb R$, and $\gr$ so that they coincide on $\kr$. This will allow us to refer to all three by the same symbol $B$, without ambiguity.

Extension of scalars identifies the space of $\Bbb R$-linear maps from $\sr$ to $\fg^\Bbb R$ with the space of $\Bbb C$-linear maps from $\fs$ to $\fg$. Also, the Cartan involution on $\sr$ equals the restriction to $\sr$ of complex conjugation with respect to the compact real form $\fs\fu(2)$ in $\fs=\fs\fl(2,\Bbb C)$. This results in the following ``dictionary"  between the spaces of homomorphisms  defined in the preceeding sections and their analogues in the present  setting:
$$
\gathered
\Hom^\Bbb R(\fs,\fg) \ \rightsquigarrow \ \Hom^\Bbb R (\fs,\Bbb C \tsize\otimes_\Bbb R \fg^\Bbb R) \ \cong \ \Hom(\fs,\fg)
\\
\Hom^{\Bbb R,\theta}(\fs,\fg) \ \rightsquigarrow \ \Hom(\fs\fu(2), \fu_\Bbb R)\,,
\endgathered
\tag6.6
$$
and similarly in the case of $\Mor(\fs,\fg)$. The Lie algebra $\fs\fu(2)$ acts on these spaces, both via the action on the values and diagonally, so the decompositions (4.6-7) have obvious counterparts. Note that the evaluation map
$$
\aligned
&\Hom^{\Bbb R,\theta}(\fs,\fg)\owns \Phi_0 \, \mapsto \, \Phi_0(e) \qquad \text{corresponds to}
\\
&\qquad\qquad \Hom(\fs\fu(2), \fu_\Bbb R) \owns \Phi_0 \ \mapsto \ \tsize \Phi_0(\frac e2 - \frac f2) - i\,\Phi_0(\frac{ie}2 + \frac{if}2)
\endaligned
\tag6.7
$$
via the translation (6.6). We now let $G$ denote the identity component of $\Aut(\fg)$, and $\UR\subset G$ the compact real form determined by $\ur$. Then
$$
\GR \ \rightsquigarrow \ G\,,\qquad \KR \ \rightsquigarrow \ \UR
\tag6.8
$$
completes our dictionary: when the substitutions (6.6-8) are made, the results of the earlier sections -- in particular lemmas 2.11 and 2.19, theorems 3.7, 5.2, and 5.4, propositions 2.21, 2.25, and 4.9, corollaries 2.22, 2.23, and 5.24 -- hold in the setting of a complex Lie algebra.

We return to the case of a real semisimple Lie algebra $\gr = \kr\oplus\pr$. As additional datum, we suppose that an involutive automorphism
$$
\sigma \, : \, \gr \ \rightarrow\ \gr \qquad\qquad (\,\sigma^2\ = \ 1_\gr\,)
\tag6.9a
$$
is fixed. It induces a pseudo-Riemannian ``Cartan decomposition" of $\gr$,
$$
\gr \ = \ \hr \tsize\oplus \qr\,,\qquad [\hr,\hr]\subset\hr\,,\ \ \ [\hr,\qr]\subset\qr \,,\ \ \ [\qr,\qr]\subset\hr\,,
\tag6.9b
$$
with $\hr$ and $\qr$ denoting, respectively, the $(+1)$ and $(-1)$-eigenspaces of $\sigma$. We shall assume that the usual Cartan involution $\theta$ preserves this decomposition -- equivalently,
$$
\theta \circ \sigma \ = \ \sigma \circ \theta \,.
\tag6.10
$$
When this fails to be the case, it can be brought about by replacing the Cartan decomposition with an appropriate $\GR$-conjugate. The involution $\sigma$ lifts to
$$
\GR \ = \ \Aut(\gr)^0\,.
\tag6.11
$$
Let $\HR\subset\GR$ be a subgroup lying somewhere between the fixed point group $G_\Bbb R^\sigma$ and its identity component,
$$
(G_\Bbb R^\sigma)^0 \ \subset \ \HR \ \subset G_\Bbb R^\sigma \,.
\tag6.12
$$
Then $\HR$ preserves the decomposition (6.9b), and thus acts on the set of nilpotents in $\qr$.

The Lie algebra $\sr = \fs\fl(2,\Bbb R)$ and its diagonal subalgebra $\ar$ furnish the simplest non-trivial example of a symmetric pair:
$$
\sigma_\fs\, : \, \sr \,\longrightarrow \, \sr\,,\qquad \sigma_\fs(\eta)\ = \ 
\cases 
\eta \ \ \ &\text{for} \ \ \eta \in \ar \\ -\eta\ \ \ &\text{for} \ \ \eta \in \Bbb R \, e \tsize\oplus \Bbb R \,f
\endcases
\tag6.13a
$$
is the involution, and
$$
\fs_\Bbb R \ =  \ \ar \, \tsize\oplus \,(\,\Bbb R \, e \oplus \Bbb R\,
f\,)
\tag6.13b
$$
the non-Riemannian Cartan decomposition; note that $\sigma_\fs$ does commute with the Cartan involution $\theta_\fs$ -- cf. (2.17a). The group $\SR = PSl(2,\Bbb R)$ and its diagonal subgroup $\AR$ play the roles of $\GR$ and $\HR$. The space $\,\Bbb R \, e \oplus \Bbb R\, f\,$ contains five nilpotent $\AR$-orbits, namely $\Bbb R_{>0}\,e\,$, $\Bbb R_{<0}\,e\,$, $\Bbb R_{>0}\,f$, $\Bbb R_{<0}\,f$, and $\{0\}$.

In the present setting, the roles of the set $\Mor^{\Bbb R,\theta}(\fs,\fg)$ and of the vector space $\Hom^{\Bbb R,\theta}(\fs,\fg)$ are played by
$$
\aligned
&\Mor^{\Bbb R,\theta,\sigma}(\fs,\fg)\ = \ \{\,\Phi_0 \in \Mor^{\Bbb R,\theta}(\fs,\fg) \mid \sigma \circ \Phi_0 = \Phi_0 \circ \sigma_\fs \,\}\,,
\\
&\Hom^{\Bbb R,\theta,\sigma}(\fs,\fg)\ = \ \{\,\Phi \in \Hom^{\Bbb R,\theta}(\fs,\fg) \mid \sigma \circ \Phi = \Phi\circ \sigma_\fs
\,\}\,.
\endaligned
\tag6.14
$$
Note that the decompositions (4.6-7) induce decompositions of $\Hom^{\Bbb R,\theta,\sigma}(\fs,\fg)$, because $\sigma_\fs$ acts trivially on the Casimir operator of $\fs$, and because $\sigma$ preserves the $\Phi_0$-image of $\fs$ in $\fg$. When $\Phi_0$ lies in $\Mor^{\Bbb R,\theta,\sigma}(\fs,\fg)$, we write $\Cal M^\sigma(\Phi_0)$ for the set of all those functions $\Phi(t)$ in $\Cal M(\Phi_0)$ which take values in $\Hom^{\Bbb R,\theta,\sigma}(\fs,\fg)$.

In the following, $\Cal O_\qr$ will denote a nilpotent $\HR$-orbit in $\qr$; there are only finitely many such orbits, and they are invariant under scaling by positive scale factors \cite{Se}. The $\GR$-translates of $\Cal O_\qr$ sweep out a nilpotent  $\GR$-orbit $\Cal O \subset \gr$. We use $\Cal O$ to normalize the Killing form, as in (2.9-10). The moment map (2.3) restricts to an $(\HR \cap \KR)$-equivariant map
$$
m \, : \, \Cal O_\qr \ \longrightarrow \ \hr \cap \pr\,,
\tag6.15
$$
the moment map associated to the $\HR$-action on $\Cal O_\qr$. As in the absolute case, the multiplicative group $\Bbb R^+$ acts, by scaling, on $\Cal O_\qr$ and on the set of critical points of the function $\zeta\mapsto \|m(\zeta)\|^2$; these actions commute with those of $\HR\cap\KR$. By definition,
$$
C(\Cal O_\qr)\ = \ \{\,\zeta \in \Cal O_\qr \ \, \mid \ \, \text{$\zeta$ is a critical point for $\|m\|^2$, and $\|\zeta\|=1$\,}\}
\tag6.16
$$
is the core of the nilpotent orbit $\Cal O_\qr$.

\proclaim{6.17 Proposition} The function $\,\zeta \mapsto \|m(\zeta)\|^2\,$ on $\,\Bbb S(\Cal O_\qr)\,$ is Bott-Morse. Its set of critical points coincides with the set minima, and consists of exactly one $(\HR\cap\KR)$-orbit. In particular, the core $C(\Cal O_\qr)$ is non-empty, and $\HR\cap\KR$ acts transitively on it. The map $\Phi \mapsto \Phi(e)$ establishes a $\HR\cap\KR$-equivariant bijection
$$
\{\, \Phi_0 \in \Mor^{\Bbb R,\theta,\sigma}(\fs,\fg) \mid \Phi_0(e) \in \Cal O_\qr \,\} \ \cong \ C(\Cal O_\qr)\,,
$$
and, at the point $\zeta = \Phi_0(e)\in C(\Cal O_\qr)$, identifies
$$
\tsize \bigoplus_{\ell\geq 2}\Hom^{\Bbb R,\theta,\sigma}(\fs,\fg(\ell))(\ell -2)
$$
$(\HR\cap\KR)_\zeta$-equivariantly with $(T_{C(\Cal O_\qr)}\Cal O_\qr)_\zeta$, the fiber at $\zeta$ of  the normal bundle of the orbit $\Cal O_\qr$ along its core. There exists a real analytic, $(\HR\cap\KR)$-equivariant isomorphism $\Cal O_\qr \simeq T_{C(\Cal O_\qr)}\Cal O_\qr\,$. Lastly, the Kronheimer diffeomorphism (3.7) induces
$$
\Cal M^\sigma(C(\Cal O_\qr))\ \cong \ \Cal O_\qr\,;
$$
here $\Cal M^\sigma(C(\Cal O_\qr))$ refers to the union of the $\Cal M^\sigma(\Phi_0)$ parameterized by those morphisms $\Phi_0$ which correspond to points in $C(\Cal O_\qr)$.
\endproclaim

These statements are analogous to (2.11), (2.19), (2.21), (4.9), (2.25), and (3.7), but not all them can be deduced directly from those results in the  absolute case. We will indicate briefly how to modify the earlier arguments so that they apply in the present situation. As noted in the proof of proposition 2.21, the fact that $\|m\|^2$ is Bott-Morse on $\,\Bbb S(\Cal O_\qr)$ is a general property of moment maps for linear actions. By (2.20) and the fact that the moment map $m$ in (6.15) is the restriction of the moment map (2.3) from $\Cal O$ to $\Cal O_\qr$, we conclude that the critical set of $\|m\|^2 : \Cal O_\qr \to \BR_{>0}$ is the intersection of $\Cal O_\qr$ with the critical set of $\|m\|^2$ in $\Cal O$. It follows that the core $C(\Cal O_\qr)$ is the intersection of the core $C(\Cal O)$ with $\Cal O_\qr$, and $C(\Cal O_\qr)$ consists precisely of the minima of $\|m\|^2$ on $\,\Bbb S(\Cal O_\qr)$. Any $\zeta\in C(\Cal O_\qr)$ can be embedded in a strictly normal $S$-triple, and by \cite{Ma} or \cite{Se, lemmas 1.4,1.5} such strictly normal $S$-triples constitute an $H_\Bbb R\cap\KR$-orbit. This proves the analogues of (2.11) and (2.21). The proof of lemma 2.19 can now be adapted to establish the bijection $\{\, \Phi_0 \in \Mor^{\Bbb R,\theta,\sigma}(\fs,\fg) \mid \Phi_0(e) \in \Cal O_\qr \,\} \ \cong \ C(\Cal O_\qr)$.

Let us explain next how to modify the proof of proposition 4.9 in the present setting. We denote the complexifications of $\hr$ and $\qr$ by $\fh$ and  $\fq$, respectively. The evaluation map $\Phi\to \Phi(e)=\Phi(\zeta)$ sends
$$
\tsize\bigoplus_{\ell\geq 2}\Hom^{\theta,\sigma}(\fs,\fg(\ell))(\ell -2) \ \longrightarrow \  \fq\,.
\tag6.18
$$
By (4.10a) the image of (6.18) lies in $[\zeta,\fg(\ell)]$ and, because $[\fq,\zeta]\subset\fh$, in $[\zeta,\fq]$. This proves the analogue
$$
\{\,\Phi(e) \mid \Phi\in \Hom^{\theta,\sigma}(\fs,\fg(\ell))(\ell-2)\,\} \ \subset \ [\zeta, \fq\cap\fg(\ell)]
\tag6.19
$$
of (4.10a). It remains to prove the analogue of (4.10b):
$$
\gathered
\text{for each $\xi\in\fh\cap\fg(\ell)$, there exist $\Phi\in \Hom^{\theta,\sigma}(\fs,\fg(\ell))(\ell-2)$ and $\eta\in\fh\cap\fk$}
\\
\text{so that}\ \Phi(e) = [e,\xi\!+\!\eta]\,; \ \text{in this situation,}\ [e,\xi] \ \text{uniquely determines} \ \Phi\,.
\endgathered
\tag6.20
$$
Statement (4.10b) implies the existence of such a $\Phi\in\Hom^{\theta}(\fs,\fg(\ell))(\ell-2)$, an $\eta\in\fk$, and the fact that $[e,\xi]$ uniquely determines $\Phi$. From (4.17b) we conclude that $\eta\in\fh\cap\fk$. Using the defining property $\Phi(e) = [e,\xi\!+\!\eta]$ and formula (4.17b) for $\Phi(h)$ one checks readily that $\Phi\in \Hom^{\theta,\sigma}(\fs,\fg(\ell))(\ell-2)$.

With the appropriate changes in notation, the proof of proposition 2.25 carries over almost word-for-word, giving the isomorphism $\Cal O_\qr \simeq T_{C(\Cal O_\qr)}\Cal O_\qr\,$.  

Vergne \cite{Ve} observed that Kronheimer's isomorphism (3.7) restricts to an isomorphism $\,\Cal M^\sigma(C(\Cal O_{\fq_\BR}))\simeq \Cal O_{\fq_\BR}$. Indeed, by definition, any $\Phi(\,\cdot\,)\in \Cal M^\sigma(C(\Cal O_\qr))$ intertwines $\sigma$ and $\sigma_\fs$, hence $\Phi(0)(e)\in\Cal O\cap\qr\,$ -- cf. (6.13b); here $\,\Cal O$ again denotes the $\GR$-orbit containing $\Cal O_{\fq_\BR}$. Since $\Cal O_\qr$ is a union of connected components of $\Cal O\cap\qr$, a continuity argument shows that Kronheimer's isomorphism (3.7) restricts to a one-to-one map
$$
\Cal M^\sigma(C(\Cal O_\qr))\ \longrightarrow \Cal O_\qr\,.
\tag6.21
$$
To see that it is onto, we consider a particular $\zeta\in\Cal O_\qr$ and the corresponding $\Phi(\,\cdot\,)\in \Cal M(C(\Cal O))$. The function $t \mapsto \sigma\circ\Phi(t)\circ \sigma_\fs$ also satisfies the defining conditions 3.6, and $\sigma\circ\Phi(0)\circ \sigma_\fs(e)= -\sigma\circ\Phi(t)(e) = \zeta$. By uniqueness, $\Phi(\,\cdot\,)= \sigma\circ\Phi(\,\cdot\,)\circ \sigma_\fs$, hence $\Phi(\,\cdot\,)\in \Cal M^\sigma(C(\Cal O_{\fq_\fs}))$. Thus (6.21) is surjective, as was to be shown.

\vskip 3\jot

\subheading{\bf 7. The Sekiguchi correspondence}
\vskip 1\jot

The Sekiguchi correspondence in its most general form establishes a bijection between nilpotent orbits attached to certain pairs of commuting involutions \cite{Se}. The complete statement and its specialization to the case of interest to us involve substantial notational overhead. For this reason, we discuss only the most important particular case; however, all statements and arguments extend readily\footnote{Except for the orientation statements in theorem 7.20, which needs to be modified when there are no complex and symplectic structures to orient the orbits in question.} to the general case.

We use the notation and conventions of the previous section. In particular, $\fg$ arises as the complexification of $\gr = \kr \oplus \pr$. We let $\theta$ denote both the Cartan involution of $\gr$ and its extension to $\fg$, and $\tau$ complex conjugation with respect to the compact real form $\ur\subset\fg$. Then, by construction of
$\ur$,
$$
\tau \zeta \ = \ \theta\bar\zeta\ = \ \overline{\theta\zeta}\qquad\qquad (\,\zeta\in\fg\,).
\tag7.1
$$
Here $\bar\zeta$ refers to the complex conjugate, relative to $\gr$. The complexification $\fk$ of $\kr$ corresponds to a subgroup
$$
K \ \subset \ G \ = \ \text{identity component of}\ \Aut(\fg)\,.
\tag7.2
$$
The complex group $G$ also contains
$$
\GR \ = \ (\Aut(\gr))^0\ \ \ \text{and}\ \ \ \UR \ = \ (\Aut(\ur))^0\,,
\tag7.3
$$
as noncompact and compact real form, respectively.

We shall consider nilpotent $K$-orbits in $\fp = \bc\otimes_\br\pr$ on the one hand, and nilpotent $\GR$-orbits in $i\,\gr$ on the other; $\Cal O_\fp$  will be the generic symbol for the former, and $\Cal O_\br$ for the latter. To avoid trivial exceptions, we always exclude the orbit $\{0\}$. Recall the definition (2.17b) of the basis $\{e,f,h\}$
of $\fs$.

In the discussion in \S2, we can make the trivial substitution of $i\,\gr$ for $\gr$. Then, as is shown there, $\,\Phi_0 \mapsto \Phi_0 (i\,e)\,$ induces a $\KR$-equivariant bijection
$$
\{\,\Phi_0 \in \Mor^{\br,\theta}(\fs,\fg) \mid \Phi_0(i\,e) \in \Cal O_\br \,\} \ \cong \ C(\Cal O_\br)\,,
\tag7.4a
$$
for every nilpotent $\GR$-orbit $\Cal O_\br\neq \{0\}$ in $\,i\,\gr$. When we look at all nilpotent $\GR$-orbits in $i\gr$ simultaneously, (7.4a) sets up a natural bijection
$$
\{\,\text{nilpotent $\GR$-orbits in $i\gr$}\,\}\ @>{\ \sim \ }>>\ \{\,\text{$\KR$-orbits in $\,\Mor^{\br,\theta}(\fs,\fg)$}\,\}
\tag7.4b
$$
We shall argue shortly that the results in \S 6 imply an analogous statement for nilpotent $K$-orbits in $\fp\,$: when $\Cal O_\fp\neq \{0\}$ is a nilpotent $K$-orbit in $\fp$, the assignment $\,\Phi_0 \mapsto \Phi_0 (\frac h2 + \frac {ie}2 + \frac {if}2)\,$ induces a $\KR$-equivariant bijection
$$
\{\,\Phi_0 \in \Mor^{\br,\theta}(\fs,\fg) \mid \Phi_0 (h + ie +if) \in \Cal O_\fp \,\} \ \cong \ C(\Cal O_\fp)\,,
\tag7.5a
$$
which, in turn, determines a bijection
$$
\{\,\text{nilpotent $K$-orbits in $\,\fp$\,}\}\ @>{\ \sim \ }>>\ \{\,\text{$\KR$-orbits in $\,\Mor^{\br,\theta}(\fs,\fg)$}\,\}\,,
\tag7.5b
$$
in complete analogy to (7.4b). Combining (7.5b) with the inverse of (7.4b), we obtain the {\it Sekiguchi correspondence}
$$
\{\,\text{nilpotent $K$-orbits in $\,\fp$\,}\}\ @>{\ \sim \ }>>\ \{\,\text{nilpotent $\GR$-orbits in $i\gr$}\,\}
\tag7.6
$$
\cite{Se}, which relates the $K$-orbit $\Cal O_\fp$ to the $\GR$-orbit $\Cal O_\br$ precisely when the inverse images of $C(\Cal O_\fp)$ and $C(\Cal O_\br)$ in $\Mor^{\br,\theta}(\fs,\fg)$ coincide.

We still need to establish (7.5a). For this purpose, we regard $(\fg,\fk)$ as symmetric pair over $\br$, with involution $\theta$. We appeal to proposition 6.17, which needs to be translated into the present setting by means of the ``dictionary" (6.6-8). To begin with,
$$
\aligned
&\Mor^{\br,\theta}(\fs,\fg) \ \ \rightsquigarrow \ \ \Mor(\fs\fu(2),\ur) \
\cong
\\
&\qquad\qquad \cong \ \Mor^\tau(\fs,\fg)\ =_{\text{def}} \ \{\,\Phi_0 \in \Mor(\fs,\fg) \mid \tau \circ \Phi_0 = \Phi_0 \circ \tau_\fs \,\}\,,
\endaligned
\tag7.7
$$
as in (6.6); here $\tau_\fs : \fs \to \fs$ stands for complex conjugation with respect to $\fs\fu(2)$. By the same dictionary,
$$
\Mor^{\br,\theta,\sigma}(\fs,\fg)\ \rightsquigarrow \ \{\,\Phi_0 \in \Mor^\tau(\fs,\fg) \mid \theta \circ \Phi_0 = \Phi_0 \circ \sigma_\fs \,\}\,,
\tag7.8
$$
since $\theta : \fg \to \fg$ now plays the role of the involution $\sigma$. A short calculation in $Sl(2,\bc)$ gives
$$
\sigma_\fs \ = \ \Ad c \circ \theta_\fs \circ \Ad c^{-1} \,,\qquad \text{with}\ \ \ c\ = \ \frac 1{\sqrt 2} \pmatrix 1 &i \\ i &1 \endpmatrix\,,
\tag7.9
$$
and $\Ad c$ commutes with $\tau_\fs$, so (7.7) is equivalent to the assignment
$$
\aligned
\Mor^{\br,\theta,\sigma}(\fs,\fg)\ \rightsquigarrow \ &\{\,\Phi_0 \in \Mor^\tau(\fs,\fg) \mid \theta \circ \Phi_0 \circ \Ad c = \Phi_0 \circ \Ad c \circ \theta_\fs \,\}
\\
= \ &\{\,\Phi_0\circ \Ad c^{-1}  \mid  \Phi_0 \in \Mor^{\theta,\tau}(\fs,\fg) \,\}\,.
\endaligned
\tag7.10
$$
These morphisms get evaluated on $e$, as in proposition 6.17. But $\,\Mor^{\theta,\tau}(\fs,\fg) = \Mor^{\br,\theta}(\fs,\fg)\,$ by (7.1), and $\Ad c^{-1}(e) = i(\frac h2 - \frac {ie}2 -\frac {if}2)$. This gives the correspondence (7.5), with $i(\frac h2 - \frac {ie}2 -\frac {if}2)$ in place of $\,\frac h2 + \frac {ie}2 + \frac {if}2\,$. Note that nilpotent $K$-orbits in $\fp$ are invariant under scaling by nonzero complex numbers -- this is clear in the case of $\gr=\fs\fl(2,\br)$, and follows in the general case by what has been said so far. Since $i$ has absolute value 1, it maps the core of an orbit to itself, so we can drop the factor $i$. Finally, complex conjugation permutes the nilpotent $K$-orbits \footnote{There are two equally natural definitions of the Sekiguchi correspondence. They are related by complex conjugation on the side of $K$-orbits, or alternatively, by multiplication by $-1$ on the side of $\GR$-orbits. Our choice of the correspondence is dictated by the application in \cite{SV2}.} in $\fp$, and this allows us to replace $\,\frac h2 - \frac {ie}2 -\frac {if}2\,$ by its complex conjugate.

When the Sekiguchi correspondence relates $\Cal O_\fp$ to $\Cal O_\BR$, various objects attached to the two orbits are naturally isomorphic -- the cores because of (7.4a) and (7.5a):
$$
C(\Cal O_\fp)\, \cong \, C(\Cal O_\BR)\,, \ \ \tsize\Phi_0(\frac h2 + \frac {ie}2 + \frac {if}2) \leftrightarrow \Phi_0(ie)\qquad (\,\Phi_0 \in \Mor^{\BR,\theta}(\fs,\fg)\,)\,.
\tag7.11
$$
This isomorphism is $\KR$-equivariant by construction, so the isotropy subgroups of $\KR$ at $\Phi_0(\frac h2 + \frac {ie}2 + \frac {if}2)$ and $\Phi_0(ie)$ coincide\footnote{This can be seen directly: if $k\in\KR$ fixes $\Phi_0(\frac h2 + \frac {ie}2 + \frac {if}2)$, it fixes the real and imaginary parts separately, which together generate $\Phi_0(\fs)$; similarly, if $k$ fixes $\Phi_0(ie)$, it fixes the $\Phi_0$-image of $\theta(ie)=-if$ and hence also the image of $[e,f]=h$.}. Proposition 4.9 identifies the normal space to $C(\Cal O_\BR)$ at the point $\Phi_0(ie)\in C(\Cal O_\BR)$, 
$$
\tsize\bigoplus_{r\geq 2} \Mor^{\br,\theta}(\fs,\fg(r))(r-2)\ \simeq \ (T_{C(\Cal O_\BR)}\Cal O_\BR)_{\Phi_0(ie)}\,,\qquad\Phi \mapsto \Phi(ie)\,,
\tag7.12
$$
equivariantly with respect to the isotropy subgroup of $\KR$ at $\Phi_0(ie)$. The analogue of (4.9) in the symmetric pair case -- which is part of proposition 6.17 -- identifies the normal space to $C(\Cal O_\fp)$ at $\Phi_0(\frac h2 + \frac {ie}2 + \frac {if}2) \in C(\Cal O_\fp)$,
$$
\gathered
\tsize\bigoplus_{r\geq 2} \Mor^{\br,\theta}(\fs,\fg(r))(r-2)\ \simeq \ (T_{C(\Cal O_\fp)}\Cal O_\fp)_{\Phi_0(\frac h2 + \frac {ie}2 + \frac {if}2)}\,,
\\
\tsize\Phi \mapsto \Phi(\frac h2 + \frac {ie}2 + \frac {if}2)\,,
\endgathered
\tag7.13
$$
again equivariantly with respect to the isotropy subgroup of $\KR$. The preceding statement involves the same ``translation" that we just used to establish (7.5a). Because of (7.12-13), the fiber of the normal bundle $T_{C(\Cal O_\fp)}\Cal O_\fp$ at $\Phi(\frac h2 + \frac {ie}2 + \frac {if}2)$ is isomorphic to the fiber of $T_{C(\Cal O_\BR)}\Cal O_\BR$ at $\Phi(ie)$ -- isomorphic as representation space for the common isotropy group. Thus (7.11) extends to a real analytic isomorphism of $\KR$-equivariant vector bundles
$$
T_{C(\Cal O_\fp)}\Cal O_\fp \ \simeq \ T_{C(\Cal O_\BR)}\Cal O_\BR\,.
\tag7.14
$$
Appealing to proposition 2.25 and its analogue for symmetric pairs, as stated in (6.17), we obtain
$$
\Cal O_\fp \ \ \simeq \ \ \Cal O_\gr\,,
\tag7.15
$$
a real analytic, $\KR$-equivariant isomorphism between the two orbits. 

Vergne \cite{Ve} deduces the existence of a $\KR$-equivariant diffeomorphism $\Cal O_\fp\simeq\Cal O_\gr$ from Kronheimer's description of nilpotent orbits, as follows. According to (3.7),
$$
\Cal M(C(\Cal O_\BR)) \ @>{\ \sim \ }>> \  \Cal O_\BR\,, \qquad \Phi(\,\cdot\,)\ \mapsto \ \Phi(0)(ie)\,,
\tag7.16
$$
is a $\KR$-equivariant diffeomorphism. The analogous statement for symmetric pairs in (6.17), translated as in (7.7-10), gives the $\KR$-equivariant diffeomorphism
$$
\Cal M(C(\Cal O_\fp)) \ @>{\ \sim \ }>> \  \Cal O_\fp\,, \qquad \tsize\Phi(\,\cdot\,)\ \mapsto \ \Phi(0)(\frac h2 + \frac {ie}2 + \frac {if}2)\,.
\tag7.17
$$
Since $\Cal O_\Bbb R$ and $\Cal O_\fp$ are related by the Sekiguchi correspondence, (7.11) implies
$$
\Cal M(C(\Cal O_\fp))\ =\ \Cal M(C(\Cal O_\Bbb R))\,.
\tag7.18
$$
The composition of (7.16-18) gives Vergne's interpretation of the Sekiguchi corres\-pondence. 

In the proof of the Barbasch-Vogan conjecture in \cite{SV2}, we were lead to a quite different geometric description of the correspondence. We fix a nilpotent $G$-orbit $\Cal O$ in $\fg$. Via the isomorphism $\fg \simeq \fg^*$ induced by the Killing form, $\Cal O$ can be viewed as a complex coadjoint orbit. As such, it carries a holomorphic symplectic form $\sigma_{\Cal O}$; in particular, $\Cal O$ has even complex dimension $2k$. The intersection $\Cal O\cap \fp$ decomposes into a union of finitely many $K$-orbits, all Lagrangian with respect to $\sigma_{\Cal O}$, hence of complex dimension $k$ \cite{KR}. Analogously, $\Cal O\cap i\gr$ is a union of finitely many $\GR$-orbits, Lagrangian with respect the real symplectic form $\operatorname{Re}\sigma_{\Cal O}$, hence of real dimension $2k$. We enumerate the two types of orbits as
$$
\Cal O\cap \fp \ =\ \Cal O_{\fp,1}\cup\dots\cup \Cal O_{\fp,d} \ \ \ \text{and} \ \ \ \Cal O\cap i\gr \ =\ \Cal O_{\Bbb R,1}\cup\dots\cup \Cal O_{\Bbb R,d}\,.
\tag7.19
$$
The complex structure orients the orbits $\Cal O_{\fp,j}$, which gives meaning to the $[\Cal O_{\fp,j}]$ as top dimensional cycles, with infinite support, in $\Cal O\cap\fp$. We had remarked already that $\sigma_\Cal O$ restricts to a purely imaginary form on $\Cal O\cap i\gr$. Thus $\frac 1 i\sigma_{\Cal O}$ defines a symplectic form on the $\GR$-orbits $\Cal O_{\Bbb R,j}$ -- one can check that this gives the same symplectic structure as the identification of $\Cal O_{\Bbb R,j}$ with a real coadjoint orbit via division by $i$ and the isomorphism $\gr \simeq \gr^*$ induced by the Killing form. We use the symplectic structure to orient the $\Cal O_{\gr,j}$, and to regard them as cycles $[\Cal O_{\gr,j}]$ in $\Cal O\cap i\gr$. We let $\ohf{*}(\dots, \Bbb Z)$ denote homology with possibly ``infinite supports" (Borel-Moore homology). Then, in view of (7.19),
$$
\aligned
\ohf{2k}(\Cal O \cap \frak p, \Bbb Z) \ &= \ \{\,\sum n_j\,[\Cal O_{\fp,j}]\ \mid \ n_j\in\Bbb Z\,\}\,,
\\
\ohf{2k}(\Cal O \cap i\gr, \Bbb Z) \ &= \ \{\,\sum n_j\,[\Cal O_{\gr,j}] \ \mid \ n_j\in\Bbb Z\,\}\,.
\endaligned
\tag7.20
$$
There are no relations among the $[\Cal O_{\fp,j}]$, respectively $[\Cal O_{\BR,j}]$, since we are dealing with top dimensional homology. This allows us to view the Sekiguchi correspondence as a specific isomorphism between the two homology groups.

Our description of the Sekiguchi correspondence amounts to a geometric passage between the two homology groups in (7.20). We define a real analytic family of diffeomorphisms
$$
f_t:\Cal O\ \longrightarrow \ \Cal O\,, \qquad  f_t(\zeta) \ = \ \operatorname{Ad}(\exp(t \operatorname{Re}\,\zeta))(\zeta)\qquad (\,t\in \BR\,)\,;
\tag7.21
$$
this agrees with the definition in \cite{SV2,\S6}, except for the change of variables $s=t^{-1}$. The images $(f_t)_*[\Cal O_{\fp,j}]$, $0\leq t < \infty$, of any $[\Cal O_{\fp,j}]$ constitute a real analytic family of cycles in the complex orbit $\Cal O$. We argue in \cite{SV2} that this family of cycles has a limit for $t\to +\infty$ for a priori reasons, and that the limit is a cycle in $\Cal O \cap i\gr$. The existence of the limit may seem surprising, since $f_t$ has exponential behavior for large $t$. At the end of this section, we shall say a few words about the notion of limit of a family of cycles, and about the argument for the existence of the limit in our situation.

\proclaim{7.22 Theorem} The assignment $\,\,c \,\mapsto \,\lim_{t\to +\infty}\,(f_t)_*\,c\,\,$ induces the Sekiguchi corres\-pondence, as map from $\,\ohf{2k}(\Cal O \cap \frak p, \Bbb Z)\,$ to $\,\ohf{2k}(\Cal O \cap i\gr, \Bbb Z)$. In other words,
$$
{\lim}_{t\to +\infty}\,(f_t)_*[\Cal O_\fp]\ \ = \ \ [\Cal O_\gr]
$$ 
whenever the $K$-orbit $\Cal O_\fp$ in $\Cal O\cap\fp$ and the $\GR$-orbit $\Cal O_\Bbb R$ in $\Cal O\cap i\gr$ are related by the Sekiguchi correspondence.
\endproclaim

This theorem plays a crucial role in our proof of the Barbasch-Vogan conjecture in \cite{SV2}, where it is stated as theorem 6.3. Its proof splits naturally in two parts. One establishes the existence of the limit and reduces its computation to two geometric lemmas about nilpotent orbits \cite{SV2,\S6}. The second part consists of the proofs of the two lemmas; these proofs occupy the last section of this paper and use the tools developed in the preceding sections.

We had mentioned earlier that our description of the Sekiguchi correspondence carries over to its most general version, which relates nilpotent orbits attached to symmetric pairs defined by commuting involutions \cite{Se}. The statement and the various steps of the proof apply in the general case after minimal changes, with one exception: in the absence of complex and symplectic structures, the orbits no longer carry natural orientations and -- as far as we know -- need not even be orientable. One can deal with this problem by considering the orbits as cycles with values in the orientation sheaf; the isomorphism (7.15) identifies the orientation sheaves of any two orbits related by the Sekiguchi correspondence. When that is done, theorem 7.22 remains correct as stated.

Let us comment briefly on the meaning of the limit in theorem 7.22 -- for a more detailed discussion of limits of cycles in general, see \cite{SV1}. When we restrict the family of cycles $\{\,(f_t)_*[\Cal O_\fp]\,\}$ to some finite interval $0 \leq t \leq a\,$, we obtain a submanifold with boundary in $[0,a]\times\Cal O$, and the boundary consists of the fibers over $0$ and $a\,$; in this situation, it is natural to think of the fiber over $a$ as the limit of the family as $t$ tends to $a$ from below. What matters here is not the smoothness of the family; it suffices that the total space and the map to the parameter interval $[0,a]$ be Whitney stratifiable. In the case of real algebraic, or more generally, subanalytic families of cycles, Whitney stratifiability is automatic. The family $\{\,(f_t)_*[\Cal O_\fp]\,\}$ fails to be subanalytic at $t=+\infty$. It does, however, belong to one of the {\it analytic-geometric categories\/} constructed by van den Dries-Miller \cite{DM}, using recent work in model theory \cite{W, DMM}. These analytic-geometric categories generalize the notion of subanalyticity, and share most of the important properties of the subanalytic category, such as Whitney stratifiability. This implies the existence of the limit in theorem 7.22; in effect, one can argue as if the family were subanalytic even at infinity. By definition, the limit cycle is supported on $F_\infty$, the fiber over $\{+\infty\}$ of the closure of $\,\{\,(t,f_t(\zeta)) \mid 0<t<\infty, \zeta \in \Cal O_\fp\,\}\,$ in $\,[0,+\infty]\times\Cal O\,$. A fairly simple argument identifies $F_\infty$ as $\Cal O \cap i\gr$ \cite{SV2, \S6}. Thus, according to (7.19), the limit cycle can only be an integral linear combination of the $\Cal O_{\BR,j}$. A normal slice to $\Cal O_{\BR,j}$ in $\Cal O$, at a generic point of $\Cal O_{\BR,j}$, intersects $(f_t)_*[(\Cal O_\fp)]$, for $t$ close to $+\infty$, with an intersection multiplicity $m_j$ not depending on $t$; here ``generic" is to be taken in the sense of the analytic-geometric category to which the family of cycles belongs. Essentially by definition, the intersection multiplicity $m_j$ is the coefficient of $\Cal O_{\BR,j}$ in the limit cycle. We argue in \cite{SV2, \S6} that the multiplicity $m_j$ can be calculated even at ``non-generic" points under certain 
circumstances. This argument reduces the statement of theorem 7.22 to the second of the two lemmas in the next section; the first lemma is a crucial ingredient of the proof of the second.

\vskip 3\jot

\subheading{\bf 8. Two Lemmas}
\vskip 1\jot

We work in the setting of the complexified Lie algebra $\fg = \bc\otimes_\br\gr$, as in \S\S6-7. We keep fixed, once and for all, a nilpotent $G$-orbit $\Cal O$ in $\fg - \{0\}$. Recall the family of real analytic maps
$$
f_t \, : \, \Cal O \ \longrightarrow \ \Cal O \,,\qquad f_t(\zeta) = \Ad(\exp(t\operatorname{Re}\zeta))(\zeta)\qquad (\,t\in\br\,)\,,
\tag8.1
$$
defined in (7.21); as was remarked earlier, this agrees with the definition in \cite{SV2, \S5}, except for the change of variables $t=s^{-1}$. Note that
$$
\operatorname{Re}\left(f_t(\zeta)\right)\ = \ \operatorname{Re}\zeta
\tag8.2
$$
for all $\zeta\in\Cal O$. It follows that $\{f_t\}$ is a one parameter group of diffeomorphisms.

Because of (6.5) and (7.1), the definition (2.7) of the moment map translates into
$$
m\,:\,\fg - \{0\} \ \longrightarrow \ i\,\ur\,, \qquad m(\zeta) \ = \ -\,\frac {[\zeta,\theta\bar\zeta]}{\|\zeta\|^2}
\tag8.3
$$
in the present situation. This map is invariant under the action of the maximal compact subgroup $\UR\subset G$. We are interested in the qualitative behavior of $\|m\|^2$ along trajectories $\{f_t(\zeta) \mid t\geq 0\}$ through points $\zeta \in\Cal O\cap \fp$. Thus we consider a particular $K$-orbit $\Cal O_\fp$ in $\Cal O\cap \fp$ and a point $\zeta \in \Cal O_\fp$. With this choice of $\zeta$ kept fixed, we write
$$
\gathered
m(f_t(\zeta)) \ = \ m(t) \ = \ m_1(t)+m_2(t)+m_3(t)\,,\ \ \ \text{with}
\\
m_1(t)\in \Bbb R \cdot \operatorname{Re}\,\zeta\,, \ \ \
m_2(t)\in\pr\cap (\operatorname{Re}\,\zeta)^\perp\,, \ \ \
m_3(t)\in i \kr\,.
\endgathered
\tag8.4
$$
This can be done because $i\,\ur = i\,\kr \oplus \pr$. Our first statement is \cite{SV2, lemma 6.28}, phrased in terms of the new variable $t=s^{-1}$.

\proclaim{8.5 Lemma} For $\,\zeta\in \Cal O_\fp$ as above and $\,t\in\br$\,, \, \ $\|m_1(t)\|^2 \ + \ \|m_3(t)\|^2\ \geq\  \|m(0)\|^2$\,.
\endproclaim

Before embarking on the proof of the lemma, we state the second one. Besides the $K$-orbit $\Cal O_\fp$ in $\Cal O\cap \fp$, it involves a $\GR$-orbit $\Cal O_\br$ in $\Cal O \cap i\,\gr$, which may or may not be related to $\Cal O_\fp$ by the Sekiguchi correspondence. We fix a point $\nu \in C(\Cal O_\br)$, which can be represented as
$$
\nu \ = \ i\,\Phi_0(e)\,, \qquad \text{with}\ \ \Phi_0 \in \Mor^{\br,\theta}(\fs,\fg)\,,
\tag8.6
$$
as in (7.4). The choice of $\Phi_0$ gives meaning to the decomposition (4.5) of $\fg$. The space
$$
\fq(\nu) \ = \ \tsize\bigoplus_{r\geq 1}\bigoplus_{\ell<r} \fg(r,\ell)
\tag8.7
$$
is a linear complement to $\,\Ker\ad(\nu)=\Ker\ad(e)\,$ in $\fg$, and is defined over $\br$. Thus, for $a>0$ sufficiently small, the map
$$
\gathered
\{\,(\xi,\eta)\in \gr \times \gr \ \mid \ \xi,\eta \in \fq(\nu)\cap \gr\,,\ \|\xi\|,\|\eta\| < a\,\} \ \ \longrightarrow \ \ \Cal O 
\\
(\xi,\eta) \ \ \longmapsto \ \ \Ad(\exp i\xi \exp \eta)(\nu)
\endgathered
\tag8.8
$$
sends its domain isomorphically to an open neighborhood of $\nu$ in $\Cal O$. In particular, shrinking $a$ further if necessary, we can make
$$
N(\nu,a) \ = \ \{\, \Ad (\exp i\xi)(\nu)\ \mid \ \xi \in \fq(\nu)\cap \gr \,,\ \|\xi\| < a \,\}
\tag8.9
$$
a normal slice to $\Cal O_\br$ in $\Cal O$ at the point $\nu$ -- in other words, a submanifold of $\Cal O$ that intersects $\Cal O_\br$ at the single point $\nu$, where the intersection is transverse. We remarked in \S7 that the orbits $\Cal O_\fp\,$, $\Cal O_\br$ carry natural orientations: the former as complex manifold, the latter as coadjoint orbit via $\Cal O_\br \owns i\zeta \mapsto \zeta \in \gr \cong \gr^*$, hence as canonically symplectic manifold \footnote{The orientation conventions are spelled out in detail in \cite{SV1, \S8}}. The orientation of $\Cal O_\fp$ in turn orients the diffeomorphic image $f_t(\Cal O_\fp)$. Our next statement is a more specific version of lemma 6.29 in \cite{SV2}, again phrased in terms of the variable $t=s^{-1}$.

\proclaim{8.10 Lemma} For $a$ sufficiently small and $t$ sufficiently large in terms of $a$, the submanifolds $N(\nu,a)$ and $f_t(\Cal O_\br)$ of $\Cal O$ meet either exactly once or not at all, depending on whether or not $\Cal O_\br$ is the Sekiguchi image of $\Cal O_\fp$. In the former situation, the intersection is transverse and has intersection multiplicity $+1$ when the orientation of $N(\nu,a)$ and the sign convention for intersection multiplicities are chosen so as to make $N(\nu,a)$ meet $\Cal O_\BR$ with multiplicity $+1$ at $\nu$.
\endproclaim

\demo{Proof of 8.5} Let us record some observations about nilpotents in $\Cal O \cap \fp$; they will be used not only here but also in the proof of the second lemma. We consider an arbitrary $\zeta\in \Cal O \cap \fp$, which we express as $\zeta=\xi +i\eta$ with $\xi,\eta\in\pr$. In particular, $\ad\, \xi:\gr\to\gr$ is self-adjoint with respect to the complex extension of the inner product (2.2). Thus
$$
\gathered
\zeta \ = \ \tsize\sum \zeta_\lambda \ = \ \xi+i\sum\eta_\lambda\,,\ \ \ \text{with $\lambda$ ranging over $\Bbb R$\,, and}
\\
\eta_\lambda\in \fg_\Bbb R^\lambda \ = \ \text{$\lambda$-eigenspace of $\operatorname{ad}\xi$}\,, \qquad \zeta_\lambda = \cases \xi+i\eta_0 &\text{if $\lambda = 0$}\,,
\\
i \eta_\lambda &\text{if $\lambda \neq 0$}\,.\endcases
\endgathered
\tag8.11a
$$
The nilpotence of $\zeta$ implies $0= B(\zeta,\zeta)= B(\xi+i\eta,\xi+i\eta)=\|\xi\|^2-\|\eta\|^2+2i(\xi,\eta)$. Also, $\fg_\br^\lambda \perp \fg_\br^\mu$ unless $\mu = \lambda$, hence
$$
\|\xi\|^2 \ = \ \|\eta\|^2 \ = \ \tsize\sum \|\eta_\lambda\|^2\ = \ \frac 1 2 \, \|\zeta\|^2\,, \ \ \ (\xi,\eta)\ = \ (\xi,\eta_0)\ = \ 0\,.
\tag8.11b
$$
Both $\xi$ and $\eta$ lie in $\pr$, i.e., the $(-1)$-eigenspace of $\theta$, hence
$$
\theta \eta_\lambda \ = \ -\eta_{-\lambda}\,, \qquad \|\eta_\lambda\|\ = \ \|\eta_{-\lambda}\|\,.
\tag8.11c
$$
All this applies to the point $\zeta$ referred to in the statement of the lemma.

We calculate $m(t)=m(f_t(\zeta))$, beginning with the definition (8.3) of the moment map:
$$
\gathered
m(t)\ = \ -\,\frac {[\Ad(\exp(t\xi))\zeta, \theta \overline{(\Ad(\exp(t\xi))\zeta)}]}{\|\Ad(\exp(t\xi))\zeta\|^2}
\\
=\ \frac {[\xi+i\tsize \sum_\lambda e^{\lambda t}\eta_\lambda, \xi-i\tsize\sum_\lambda e^{\lambda t}\eta_{-\lambda}]}{\sum_\lambda e^{2\lambda t}\|\zeta_\lambda\|^2}
\\
=\ \frac {\tsize\sum_{\lambda,\mu} e^{(\lambda+\mu)t}[\eta_\lambda,\eta_{-\mu}]\ - \ i\sum_\lambda \lambda(e^{\lambda t}+e^{- \lambda t})\eta_\lambda}{\sum_\lambda e^{2\lambda t}\|
\zeta_\lambda\|^2}\,.
\endgathered
$$
The imaginary part of this expression equals $m_3(t)$, and $i\eta_\lambda= \zeta_\lambda$ if $\lambda\neq 0$. We conclude:
$$
\gathered
m_3(t) \ = \ -\, \frac { \sum_\lambda \lambda(e^{\lambda t}+e^{- \lambda t})\zeta_\lambda}{\sum_\lambda e^{2\lambda t}\|\zeta_\lambda\|^2}\,,
\\
\|m_3(t)\|^2\ = \ \frac { \sum_\lambda \lambda^2(e^{\lambda t}+e^{- \lambda t})^2\|\zeta_\lambda\|^2}{(\sum_\lambda e^{2\lambda t}\|\zeta_\lambda\|^2)^2}\,.
\endgathered
\tag8.12
$$
On the other hand,
$$
\gathered
m_1(t) \ = \ \frac{(\operatorname{Re}\,m(t), \xi)}{\|\xi\|^2}\ \xi \ = \ 2\ \frac{B(\operatorname{Re}\,m(t), \xi)}{\|\zeta\|^2}\ \xi
\\
= \ 2\ \frac{\tsize\sum_{\lambda,\mu}e^{(\lambda+\mu)t}B([\eta_\lambda,\eta_{-\mu}],
\xi)}{\|\zeta\|^2\ \sum_\lambda e^{2\lambda t}\|\zeta_\lambda\|^2}\ \xi\ = \ -2\ \frac{\tsize\sum_{\lambda,\mu}e^{(\lambda+\mu)t}B(\eta_\lambda, [\xi,\eta_{-\mu}])}{\|\zeta\|^2\ \sum_\lambda e^{2\lambda t}\|\zeta_\lambda\|^2}\ \xi 
\\
= \ 2\ \frac{\tsize\sum_{\lambda}\lambda e^{2\lambda t}B(\eta_\lambda,\eta_{-\lambda})}{\|\zeta\|^2\ \sum_\lambda e^{2\lambda t}\|\zeta_\lambda\|^2}\
\xi\ = \ 2\ \frac{\tsize\sum_{\lambda}\lambda e^{2\lambda t}\|\eta_\lambda\|^2}{\|\zeta\|^2\ \sum_\lambda
e^{2\lambda t}\|\zeta_\lambda\|^2}\ \xi\,.
\endgathered
$$
In the numerator, only the summands corresponding to non-zero $\lambda$ matter, so we can replace $i\eta_\lambda$ by $\zeta_\lambda$\,, giving us
$$
\gathered
m_1(t)\ = \ 2\ \frac{\tsize\sum_{\lambda}\lambda e^{2\lambda t}\|\zeta_\lambda\|^2}{\|\zeta\|^2\ \sum_\lambda e^{2\lambda t}\|\zeta_\lambda\|^2}\ \xi\,,
\\
\|m_1(t)\|^2\ = \ 2\ \frac{(\tsize\sum_{\lambda}\lambda e^{2\lambda t}\|\zeta_\lambda\|^2)^2}{\sum_\lambda \|\zeta_\lambda\|^2\ (\sum_\lambda e^{2\lambda t}\|\zeta_\lambda\|^2)^2}\,;
\endgathered
\tag8.13
$$
in the second line, we have used the equality $\sum_\lambda \|\zeta_\lambda\|^2 = \|\zeta\|^2= 2 \|\xi\|^2$. At $t=0$, $m(0)$ is a positive multiple of $-[\zeta,\theta\bar\zeta]$, which lies in $i\kr$. Thus
$$
\gathered
m(0) \ = \ -2\ \frac { \sum_\lambda \lambda\zeta_\lambda}{\sum_\lambda \|\zeta_\lambda\|^2}\,,
\\
\|m(0)\|^2 \ = \ 4\ \frac { \sum_\lambda \lambda^2\|\zeta_\lambda\|^2}{(\sum_\lambda \|\zeta_\lambda\|^2)^2}\,,
\endgathered
\tag8.14
$$
as follows from (8.12) with $t=0$.

To prove the inequality asserted by the lemma, we rewrite it in terms of the expressions (8.12-14) and clear the (positive) denominators. This transforms the inequality into the following equivalent form:
$$
\gathered
2(\tsize\sum_\lambda \|\zeta_\lambda\|^2)(\tsize\sum_{\lambda}\lambda e^{2\lambda t}\|\zeta_\lambda\|^2)^2 \ +\ (\sum_\lambda \| \zeta_\lambda\|^2)^2 (\sum_\lambda
\lambda^2(e^{\lambda t}+e^{- \lambda t})^2\|\zeta_\lambda\|^2)
\\
\geq\ 4(\tsize\sum_\lambda \lambda^2\|\zeta_\lambda\|^2)  (\sum_\lambda e^{2\lambda t}\|\zeta_\lambda\|^2)^2\,.
\endgathered
\tag8.15
$$
The original inequality is homogenous in $\zeta$. So is the inequality (8.15) when one allows only scaling by real numbers and gives $\lambda$ -- which is a typical eigenvalue of $\ad\,\xi= \operatorname{ad}(\operatorname{Re}\zeta)$ -- the same weight as $\|\zeta_\lambda\|$. Thus we are free to renormalize $\zeta$, subject to the condition $2=\|\zeta\|^2= \sum_\lambda \| \zeta_\lambda\|^2$. We set $a_\lambda= \|\zeta_\lambda\|^2$. Then $\sum_\lambda a_\lambda=2$, and $a_0=\|\xi\|^2 + \|\eta_0\|^2 \geq 1$ by (8.11); also, $a_{-\lambda}=a_\lambda$, again by (8.11). We note that $(e^{\lambda t}+e^{- \lambda t})^2= e^{2\lambda t} + 2 + e^{-2\lambda t}$, and replace $2t$ by $t$ throughout. At this point, the inequality to be proved becomes
$$
\gathered
(\tsize\sum_{\lambda}\lambda e^{\lambda t} a_\lambda)^2 \ +\  \sum_\lambda \lambda^2 (e^{\lambda t}+2+e^{- \lambda t})a_\lambda \geq\ (\tsize\sum_\lambda \lambda^2a_\lambda)  (\sum_\lambda e^{\lambda t}a_\lambda)^2\,,
\\
\text{subject to the conditions} \ \ a_\lambda=a_{-\lambda}>0\,, \ \ a_0\geq 1\,, \ \ \tsize \sum_\lambda a_\lambda = 2\,.
\endgathered
\tag8.16
$$
There must be at least one pair of non-zero indices $\pm\lambda$; otherwise $\xi$ and $\eta$ would commute, making $\zeta$ semisimple -- impossible, since $\zeta$ is a non-zero nilpotent. One further reformulation of the inequality to be proved: we define
$$
h(t) \ = \ \tsize\sum_\lambda \,a_\lambda\,e^{\lambda t}\,, \qquad t\in
\Bbb R\,.
\tag8.17a
$$
This transforms the inequality into the form
$$
h'(t)^2 \ +\ 2h''(t)\ + \ 2h''(0) \ \geq \ h''(0)h(t)^2\,,
\tag8.17b
$$
with the $a_\lambda$ still subject to the conditions listed in (8.16).

The function $h(t)$ has a globally convergent Taylor series. We can therefore verify (8.17b) by establishing the corresponding inequalities for all derivatives, at $t=0$, of the expressions on both sides, including the 0-th derivative, of course. The conditions on the $a_\lambda$ imply, in particular,
$$
\alignedat2
&\text{a)} \ \ h(0)  \ = \ \tsize\sum_{\lambda} \, a_\lambda\ = \ 2\,;
\\
&\text{b)} \ \ 0<\tsize\sum_{\lambda \neq 0} \, a_\lambda\ \leq 1\,;
\\
&\text{c)} \ \ h^{(2k)}(0)\  = \ \tsize\sum_{\lambda \neq 0} \, \lambda^{2k}\,a_\lambda\,, \ \ \ \ \ \ &&\text{for $k>0$}\,;
\\
&\text{d)} \ \ h^{(2k+1)}(0)\  = \ 0\,, \ \ \ &&\text{for $k\geq0$}\,.
\endalignedat
\tag8.18
$$
This gives the inequality at $t=0$, as an equality, in fact. We still must show that
$$
\frac{d^{(k)\ }}{dt^k}\left( h'(t)^2 \ +\ 2h''(t)\ + \ 2h''(0)\right)|_{t=0} \ \ \geq\ \ \ \frac{d^{(k)\ }}{dt^k}\left(h''(0)h(t)^2\right)|_{t=0}
$$
for all $k>0$, or equivalently,
$$
\aligned
&\sum_{\ell = 0}^k \binom k \ell h^{(\ell+1)}(0)h^{(k-\ell+1)}(0) \ \ + \ \ 2 h^{(k+2)}(0)
\\
&\geq\ \ \sum_{\ell = 0}^k \binom k \ell h^{(2)}(0) h^{(\ell)}(0)h^{(k-\ell)}(0)\,.
\endaligned
\tag8.19
$$
When $k$ is odd, both sides reduce to zero because of (8.18d). To deal with the even case, we replace $k$ by $2k$, omit the odd derivatives in the two sums, and separate out the summands involving $h(0)=2$.  This reduces the problem to showing that
$$
\gathered
\sum_{\ell = 1}^{k} \binom {2k}{2\ell-1} h^{(2\ell)}(0)h^{(2k-2\ell+2)}(0) \ \ + \ \ 2 h^{(2k+2)}(0)
\\
\geq\ \ \sum_{\ell = 1}^{k-1} \binom {2k} {2\ell} h^{(2)}(0) h^{(2\ell)}(0)h^{(2k-2\ell)}(0) \ \ +\ \ 4 h^{(2)}(0) h^{(2k)}(0)\,,
\endgathered
\tag8.20
$$
still for $k>0$.

We shall reorganize the terms on both sides of (8.20) and then compare corresponding terms, using the Chebychev inequality. First the left hand side of (8.20):
$$
\aligned
&\sum_{\ell = 1}^{k} \binom {2k}{2\ell-1} h^{(2\ell)}(0)h^{(2k-2\ell+2)}(0) \ \ + \ \ 2 h^{(2k+2)}(0)
\\
=\ \
&\sum_{\ell = 1}^{k} \left\{\binom {2k-1}{2\ell-1} +\binom {2k-1}{2\ell-2} \right\}h^{(2\ell)}(0)h^{(2k-2\ell+2)}(0)\ \ + \ \ 2 h^{(2k+2)}(0)
\\
=\ \ &\sum_{\ell = 1}^{k}\binom {2k-1}{2\ell-1} h^{(2\ell)}(0)h^{(2k-2\ell+2)}(0)\ \ +
\\
&\qquad \qquad+\ \sum_{\ell = 0}^{k-1} \binom {2k-1}{2\ell} h^{(2\ell+2)}(0)h^{(2k-2\ell)}(0)\ \  + \ \ 2 h^{(2k+2)}(0)
\\
=\ \
&\sum_{\ell = 1}^{k-1}\binom {2k-1}{2\ell-1} h^{(2\ell)}(0)h^{(2k-2\ell+2)}(0) \ \ +
\\
&\qquad\qquad +\ \sum_{\ell = 1}^{k-1} \binom {2k-1}{2\ell} h^{(2\ell+2)}(0)h^{(2k-2\ell)}(0)\ \  +
\\
&\qquad\qquad\qquad\qquad+ \ \  2h^{(2)}(0)h^{(2k)}(0) \ \  + \ \ 2 h^{(2k+2)}(0)\,.
\endaligned
\tag8.21
$$
Now the right hand side:
$$
\aligned
&\sum_{\ell = 1}^{k-1} \binom {2k} {2\ell} h^{(2)}(0) h^{(2\ell)}(0)h^{(2k-2\ell)}(0) \ \ +\ \ 4 h^{(2)}(0) h^{(2k)}(0)
\\
=\ \ &\sum_{\ell = 1}^{k-1} \left\{\binom {2k-1}{2\ell-1} +\binom {2k-1}{2\ell} \right\}h^{(2)}(0) h^{(2\ell)}(0)h^{(2k-2\ell)}(0) \ \ +
\\
&\qquad\qquad+\ \ 4 h^{(2)}(0) h^{(2k)}(0)
\\
=\ \ &\sum_{\ell = 1}^{k-1} \binom {2k-1}{2\ell-1} h^{(2)}(0) h^{(2\ell)}(0)h^{(2k-2\ell)}(0)\ \ +
\\
&\qquad\qquad + \ \
\sum_{\ell = 1}^{k-1} \binom {2k-1}{2\ell} h^{(2)}(0) h^{(2\ell)}(0)h^{(2k-2\ell)}(0) \ \ +
\\
&\qquad\qquad\qquad\qquad+\ \ 4 h^{(2)}(0) h^{(2k)}(0)\,.
\endaligned
\tag8.22
$$
Matching up corresponding terms on the right in  (8.21-22), we see that the inequality (8.20) can be reduced to
$$
h^{(2\ell+2)}(0) \ \ \geq \ \ h^{(2)}(0)h^{(2\ell)}(0)
\tag8.23
$$
for all $\ell$. This is equivalent to
$$
\sum _{\lambda\neq 0}\, \lambda^{2\ell+2}\, a_\lambda \ \ \geq \ \ \left(\sum_{\lambda\neq 0}\, \lambda^{2}\, a_\lambda\right)\left(\sum_{\lambda\neq 0}\, \lambda^{2\ell}\, a_\lambda\right)\,,
\tag8.24
$$
because of (8.18c). We now appeal to Chebychev's inequality as stated in \cite{HLP, (2.17.1)}, for example:
$$
\frac{\sum _{\lambda\neq 0}\, \lambda^{2\ell+2}\, a_\lambda} {\sum_{\lambda\neq 0}\,a_\lambda}\ \  \geq \ \ \left(\frac{\sum _{\lambda\neq 0}\, \lambda^{2}\, a_\lambda} {\sum _{\lambda\neq 0}\,a_\lambda}\right) \left(\frac{\sum_{\lambda\neq 0}\, \lambda^{2\ell}\, a_\lambda} {\sum _{\lambda\neq 0}\,a_\lambda}\right)\,.
\tag8.25
$$
But $\,0<\sum _{\lambda\neq 0}a_\lambda\leq1$ by (8.18b), so (8.25) implies (8.24),
and hence lemma 8.5.
\enddemo

\demo{Proof of 8.10} We express the point $\nu$ as in (8.6) and use the morphism $\Phi_0$ to identify $\fs=\fs\fl(2,\Bbb C)$ with a $\theta$-stable, conjugation invariant subalgebra of $\fg$. In particular,
$$
\nu \ = \ i\,e\,.
\tag8.26
$$
We must show: for $a>0$ sufficiently small and $t>0$ sufficiently large, the equation
$$
f_t(\zeta) \ = \ \Ad(\exp i\kappa)(i\,e)\,,\qquad \text{with}\ \ \zeta \in \Cal O_\fp\,,\ \ \kappa \in \fq(ie)\cap\gr\,,\ \ \|\kappa\| < a\,,
\tag8.27
$$
has exactly one solution when $\Cal O_\fp$ and $\Cal O_\BR$ are Sekiguchi related, and no solution otherwise.

It is easy to produce a solution when it is supposed to exist. Thus, for the moment, we assume that the two orbits {\it are\/} related. Note that the identity $se^{2st}=1$, with $t \geq 0\,,\ 0 < s \leq 1\,$, implicitly describes $s=s(t)$ as a decreasing function of $t$, and $\lim_{t\to\infty}s(t)=0$. A simple calculation in $SL(2,\BC)$ shows:
$$
\aligned
f_t(s\,h \,+\, i\,s\,e \,+\, i\,s\,f)\ &= \ \Ad(\exp(s\,t\,h))(s\,h \,+ \, i\,s\,e \,+\, i\,s\,f) 
\\
&= \ s\,h \,+\, i\,e \,+\, i\,s^2f \ = \ \Ad(\exp (i\,s\,f))(i\,e)\,,
\\
\text{and}\,\ s(h\,+\,i\,e\,+\,i\,f)\ \, &\text{lies in the $K$-orbit related to the $\GR$-orbit of}\ \, i\,e\,;
\endaligned
\tag8.28a
$$
in other words, the relation (8.27) with $\zeta = s(h+ie+if)$ and $\kappa = sf$ -- which does lie in $\fq(ie)\cap \gr$, as required. With little more effort, one checks that
$$
\gathered
\text{in the case of $(\gr,\kr)= (\fs\fl(2,\BR),\fs\fo(2))$, with $t>0$,}
\\
\text{the above solution of the equation (8.27) is the only solution}
\\
\text{with the property that $\,\xi \in \BR\,h\,$ and $\,\kappa \in \BR\,f\,$.}
\endgathered
\tag8.28b
$$
In fact, for $(\gr,\kr)= (\fs\fl(2,\BR),\fs\fo(2))$ and $t>0$, it is the only solution, even without the additional hypotheses on $\xi$ and $\kappa$, as will follow from the arguments below. We shall need to know certain properties of the solution (8.28a):
$$
\gathered
m(s\,h \,+\, i\,e \,+\, i\,s^2f)\ = \ (1+s^2)^{-1}(\,(1-s^2)\,h \,-\, 2\,i\,s\,e \,+\, 2\,i\,s\,f)\,,
\\
\|s\,h \,+\, i\,e \,+\, i\,s^2f\| \ = \ 1+s^2\,,\qquad \|m(s\,h \,+\, i\,e \,+\, i\,s^2f)\|^2\ = \ 2\,;
\endgathered
\tag8.29
$$
this follows from the description (8.3) of the moment map and another simple calculation.

In the general situation, let us suppose that (8.27) does have a solution, with $a>0$ sufficiently small and $t>0$ -- the meaning of ``sufficiently small" will be specified later. We write $\zeta = \xi + i\eta$, as in (8.11), and we define
$$
s \ = \ \frac {\|\xi\|}{\sqrt 2}\,.
\tag8.30
$$
The present meaning of $s$ appears to be different from that in (8.28a); after the fact, we shall see that they agree. Inductively, we shall produce bounds
$$
\|\xi - sh\|\ < \ C^k s^k\,,\ \ \ \|\kappa - sf\|\ < \ C^k s^k\,,
\tag8.31
$$
for all $k\geq 1$, with some positive constant $C$ which is independent of both $k$ and $t$. For $a$ small and $\|\kappa\|<a$,
$$
\|\xi\|\ = \ \|\operatorname{Re}f_t(\zeta)\|\ \leq \ \|f_t(\zeta)\,-\,i\,e\| \ = \ \|\Ad(\exp i\kappa)(i\,e)\, - \, i\,e\|
\tag8.32
$$
is small as well. Thus we can force $\,Cs<1$, in which case (8.31) implies $\xi = sh \in \BR\,h$ and $\kappa = sf \in \BR\,f$, hence $\xi\,,\,\kappa \in \fs$. But then $\fs$ also contains $\zeta = f_{-t}(\Ad\exp(i\kappa)(ie))$; recall: $\{f_t\}$ is a one-parameter group of diffeomorphisms. Because of (8.28b), our hypothetical solution must coincide with the solution (8.28a) -- in particular, no solution exists unless the two orbits are Sekiguchi related.

At this point, we still need to establish the bounds (8.31) and to pin down the nature of the intersection of $f_t(\Cal O_\BR)$ with the normal slice -- transverse, with sign $+1$. The latter is a separate matter, and we shall deal with it last. 

To prepare for the verification of (8.31), we re-write the right hand side of (8.27). Since $\,\ad\, f(e) = -h$ and $(\ad\, f)^2(e)= -2f$,
$$
\aligned
&\Ad (\exp i\kappa)(ie) \ = \ {\sum}_{\ell\geq 0}\,\frac {i^{\ell +1}} {\ell !} \left(s\,\ad\, f \,+\,\ad( \kappa- sf)\right)^\ell (e) 
\\
&\qquad = \ i\,e \, + \, s\,h \, + \, i\,s^2 f \, + \, [e,\kappa -sf]
\\
&\qquad \ \  - \, \frac i2 \, \left( (\ad(\kappa -sf ))^2 \, + \, \ad(\kappa -sf ) \ad(sf) \, + \, \ad(sf) \ad(\kappa -sf )\right)(e)
\\
&\qquad \ \  + \ {\sum}_{\ell>2}\,\frac {i^{\ell +1}} {\ell !} \left(s\,\ad\, f \,+\,\ad( \kappa- sf)\right)^\ell (e)\,.
\endaligned
\tag8.33a
$$
We make $a$ small enough to force $s< 1$ and $\|\kappa - sf\| < 1$. For $k>2$, $(\ad\, f)^k(e)= 0$. Thus, when we expand $\left(s\, \ad f + \ad(\kappa- sf)\right)^k \!(e)$ as a sum of monomials, every non-zero term involves at least one power of $\,\ad( \kappa- sf)$. We can therefore choose $D>0$ so that
$$
\aligned
&\| \, {\sum}_{\ell>2}\,\frac {i^{\ell +1}} {\ell !} \left(s\,\ad\, f \,+\,\ad( \kappa- sf)\right)^\ell (e)\,\| 
\\
&\qquad\qquad\qquad\qquad\qquad < \ D\,\|\kappa - sf\|\,\max (s^2,\|\kappa - sf\|^2)\,.
\endaligned 
\tag8.33b
$$
Taking the real and imaginary parts of $f_t(\zeta)= \xi + i\,\Ad(t\xi)(\eta) = \Ad (\exp i\kappa)(ie)$, we find
$$
\aligned
\text{a)}\ \ \ &\| \xi - sh - [e,\kappa -sf]\| \ \ < \ \ D\,\|\kappa - sf\|\,\max (s^2,\|\kappa - sf\|^2)\,,
\\
\text{b)}\ \ \ &\|\Ad(t\xi)(\eta) - e - s^2f\| \ \ < \ \ D\,\|\kappa - sf\|\,\max (s,\|\kappa - sf\|)\,,
\endaligned
\tag8.34
$$
now with a possibly larger value of $D$.

We remarked already that $\|\xi\|$ and $s$ are necessarily small when $a$ is small. Also, the operator $\ad\, e$ is injective on the space $\fq(ie)$, which contains both $f$ and $\kappa$. Hence $\|[e,\kappa -sf]\|$ can be bounded from below by a positive multiple of $\|\kappa -sf\|$. Using (8.34a), we now conclude:
$$
\|\xi -sh\|\ \ \ \text{and} \ \ \ \|\kappa -sf\|\ \ \ \text{are mutually bounded}
\tag8.35
$$
when $a$ is sufficiently small. In particular, this makes the two inequalities in (8.31) equivalent to one another. The first holds vacuously when $k=1$, hence so does the other.

For the inductive step, we assume that (8.31) is satisfied for some $k \geq 1$. Enlarging the constant $D$ in (8.34) if necessary -- independently of $k$ -- we can arrange
$$
\|\xi \, + \, i \Ad(t\xi)(\eta) \, -\, s\,h \,-\, i\,e \,-\, i\,s^2f\|\ < \ D\|\kappa \, - \, s\,f\| \ < \ C^kDs^k\,. 
\tag8.36
$$
But $\xi +  i \Ad(t\xi)(\eta) = f_t(\zeta)$, and $(1+s^2)^{-1}(sh + ie + is^2f)$ lies in the core $C(\Cal O)$; indeed, according to (8.29), $(1+s^2)^{-1}(sh + ie + is^2f)$ has unit length, and there the function $\|m\|^2$, which is invariant under scaling of the argument, assumes the minimum value 2. Thus (8.36) implies
$$
\dist(\,(1-s^2)^{-1} f_t(\zeta) \,,\,C(\Cal O)\,)\ \leq \ \|\,\frac {f_t(\zeta)}{1+s^2}\,-\, \frac{sh +ie + is^2}{1+s^2}\,\|\ < \ C^kDs^k\,.
\tag8.37
$$
The function $\|m\|^2 : \Bbb S(\Cal O) \to \Bbb R_{>0}$ is Bott-Morse, with minimal value 2, assumed precisely on the core. Using (8.37) and the invariance of $m$ under scaling of the argument, we find
$$
\|m(f_t(\zeta))\|^2 \ - \ 2 \ < \ C^{2k} D^2 s^{2k}\,,
\tag8.38
$$
possibly after increasing $D$, again independently of $k$. On the other hand, according to lemma (8.5),
$$
\aligned
\|m(f_t(\zeta))\|^2\ \ &= \ \|m_1(t)\|^2 \ + \ \|m_2(t)\|^2 \ + \ \|m_3(t)\|^2\
\\
&\geq \ \|m_1(t)\|^2 \ + \ \|m_3(t)\|^2 \ \geq \ \|m(\zeta)\|^2 \ \geq \ 2\,.
\endaligned
\tag8.39
$$
Combining (8.38-39), we find
$$
\|m_2(t)\| \ < \ C^kDs^k\,.
\tag8.40
$$
The moment map is differentiable, so (8.37) implies a bound on the distance between $m(f_t(\zeta))$ and $m(sh+ie+is^2f)$,
$$
\|m(f_t(\zeta)) \ - \ m(sh+ie+is^2f)\|\ < \ C^kDs^k\,,
\tag8.41
$$
with a larger $D$, if necessary. By definition of the $m_j(t)$,
$$
m_1(t) \ = \ \operatorname{Re}\,(m_1(t)\,+\,m_3(t))\ = \ \operatorname{Re}\,(m(f_t(\zeta))\,-\,m_2(t))\,.
\tag8.42
$$
At this point, we can conclude that
$$
\|\,m_1(t) \ - \ \frac{1-s^2}{1+s^2}\,h\,\| \ < \ 2C^kDs^k\,,
\tag8.43
$$
by combining the formula (8.29) for $m(sh+ie+is^2f)$ with (8.40-42).

Recall that $m_1(t)$ is a real multiple of $\operatorname{Re}\zeta = \xi$ -- a positive multiple, as follows from the explicit formula (8.13) in conjunction with (8.11):
$$
\xi \ = \ \|\xi\| \, \frac{m_1(t)}{\|m_1(t)\|}\,.
\tag8.44
$$
In this formula, we can approximate $m_1(t)$ by $(1+s^2)^{-1}(1-s^2)h$, at the expense of introducing an error term slightly larger than that in (8.43), multiplied by $\|\xi\|$. Since the inner product was normalized by the formula $\|h\|^2=B(h,h)=2$,
$$
\|\,\xi \ - \ \frac {\|\xi\|} {\sqrt 2} \, h\,\| \ < \ 3\|\xi\|C^kDs^k\,,
\tag8.45
$$
provided $s$ is sufficiently small -- which, we had seen, can be arranged by making $a$ small. We substitute $\|\xi\| = \sqrt 2 s$ -- cf. (8.30) -- and choose $C$ at least as large as $3D\sqrt 2$, giving us
$$
\|\,\xi \, - \, s\,h\,\| \ < \ C^{k+1}s^{k+1}\,.
\tag8.46
$$
In view of (8.35) this completes the inductive verification of (8.31). We had remarked already that (8.31) implies the first part of the lemma.

Now let $\Cal O_\BR$, $\Cal O_\fp$ be orbits related by the Sekiguchi correspondence, and $\nu$ a point in the core $C(\Cal O_\BR)$. We use the notation (8.26-28); in particular, we again identify $\fs$ with a subalgebra of $\fg$ and the point $\nu$ with $ie$. To shorten formulas, we set
$$
\nu_t\ = \ s(h\,+\, i\,e\, + \, i\,f) \qquad\qquad (\,0 \leq t<\infty\,)\,,
\tag8.47
$$
with $s=s(t)$ determined implicitly by $se^{2st}=1$ as before. Then $\lim_{t\to\infty} s(t) = 0$, $s(0)=1$, and
$$
f_t(\nu_t) \ = \ i\,e\, +\, s\,h \, + \,i\,s^2f\,,
\tag8.48
$$
as in (8.28). We regard tangent spaces to (real) submanifolds of $\fg$ as vector subspaces of $\fg^\BR$, i.e., of $\fg$ considered as vector space over $\BR$. However, we shall not dwell on the distinction between $\fg$ and $\fg^\BR$ from now on. We shall show:
$$
\gathered
\text{the limit of vector spaces}\ \ \ {\lim}_{t\to\infty}\ (f_t)_* \left( T_{\nu_t} \Cal O_\fp\right)
\\
\text{exists and equals}\ \ \ T_{ie}\Cal O_\BR\,.
\endgathered
\tag8.49
$$
Since $\Cal O_\BR$ and the normal slice $N(ie,a)$ meet transversely at $ie$ by construction, $f_t(\Cal O_\fp)$ must then meet $\Cal O_\BR$ transversely at $\nu_t$ for $t$ large, as asserted by the lemma.

The point $(2s)^{-1}\nu_t$ lies in the core $C(\Cal O_\fp)$. Since scaling by a positive number preserves $\Cal O_\fp$, the tangent spaces to $\Cal O_\fp$ at $\nu_t$ and $(2s)^{-1}\nu_t$ are naturally isomorphic; indeed, they are equal as subspaces of $\fg^\BR$. Appealing to (4.9) and (6.17), we find
$$
\gathered
T_{\nu_t}\Cal O_\fp \ = \ T_{\nu_t} (\KR\cdot \nu_t) \ \tsize\bigoplus \ \fd_\BR(\Phi_0)(h+ie+if) \,,
\\
\text{with}\ \ \ \fd_\BR(\Phi_0)\ = \ \tsize\bigoplus_{r\geq 2} \Hom^{\BR,\theta}(\fs,\fg(r))(r-2)\,.
\endgathered
\tag8.50
$$
We shall apply the differential of $f_t$ separately to the various summands in this decomposition of $T_{\nu_t}\Cal O_\fp$ and then take the limit as $t\to\infty$.

The map $f_t$ is $\GR$-invariant by definition. It follows that $(f_t)_*$ maps the tangent space $T_{\nu_t} (\KR\cdot \nu_t)$ isomorphically onto $T_{f_t(\nu_t)} (\KR\cdot f_t(\nu_t))$. Since $f_t(\nu_t)\to ie$, we can let $t$ tend to infinity and conclude
$$
{\lim}_{t\to\infty}\ (f_t)_* \left(  T_{\nu_t}(\KR \cdot \nu_t) \right) \ = \ T_{ie} \left( \KR \cdot ie \right)\,,
\tag8.51
$$
provided the family of $\KR$-orbits $\KR\cdot f_t(\nu_t)= \KR\cdot(ie+sh+is^2f)$ depends smoothly on $s=s(t)$ even at $s=0$. To see this, note that any $k\in\KR$ that fixes $f_t(\nu_t)$ must fix the real and imaginary parts separately, but those generate $\fs$ as Lie algebra. Similarly, if $k\in\KR$ fixes $ie$, it must fix also $if=-i\theta e$, which  together with $ie$ generates $\fs$. The constancy of the isotropy subgroups of $(\KR)_{f_t(\nu_t)}=(\KR)_\fs$ even at $s=0$ implies the smooth dependence of the $\KR$-orbits, hence (8.51).

Recall the decomposition $\fg = \oplus_{r,\ell}\,\fg(r,\ell)$ defined in \S4. For $\eta \in \fg(r)$, we let $\eta_\ell$ denote the component of $\eta$ in $\fg(r,\ell)$. We shall need to know:
$$
\aligned
\text{a)}\ \ &\text{the map}\ \ \Phi \mapsto \Phi\,h\ \ \text{establishes an isomorphism}
\\
&\Hom^{\BR,\theta}(\fs,\fg(r))(r-2)\ \ \simeq \ \ \tsize\bigoplus_{0 \leq \ell < r}\left(\fg(r,\ell)\,+\,\fg(r,-\ell) \right) \cap \pr\,;\vspace{-1.5\jot}
\\
\text{b)}\ \ &\Phi\in \Hom(\fs,\fg(r))(r-2) \qquad \Longrightarrow \qquad (\Phi\,h)_\ell \ = \ 0 \ \ \text{if}\ \ \ell=\pm r\,, \ \ \ \text{and}
\\
&(\Phi\,e)_{\ell+2}\ = \ \frac {-1}{r-\ell}\,[\,e\,,\,(\Phi\,h)_\ell\,]\,, \ \ \ (\Phi\,f)_{\ell-2}\ = \ \frac 1{r+\ell}\,[\,f\,,\,(\Phi\,h)_\ell\,]\,.
\endaligned
\tag8.52
$$
The assertion b) is established in \cite{S,\S9}, in the arguments\footnote{The hypotheses ``if $r=n$ or $r=n-2$, and if $s=\pm n,\,\pm(n-2)$" in \cite{S,(9.53)} are irrelevant in the present setting; in other words, one should argue as in \cite{S}, but with $X_n=Y_n=Z_n=0$\.} leading up to (9.53) in that paper; alternatively, one can deduce b) directly from the identity (4.14b) in the proof of proposition 4.9. Because of b), the map $\Phi\mapsto \Phi\,h$ is certainly injective on $\Hom^{\BR,\theta} (\fs,\fg(r))(r-2)$, and $\Phi\,h$ has no components in $\fg(r,\pm r)$. But any such $\Phi$ respects the Cartan decomposition and real structure, so $\Phi\,h$ lies in $\fg(r)\cap \pr$. The space $\fg(r)\cap \pr$ is invariant under $(\ad\, h)^2$, hence splits into the direct sum of the subspaces $\left(\fg(r,\ell)\,+\,\fg(r,-\ell) \right) \cap \pr\,$. Since $(\Phi\,h)_{\pm r}=0$, the map $\Phi \mapsto \Phi\,h$ in a) does take values in $\bigoplus_{0 \leq \ell r} \left(\fg(r,\ell)\,+\,\fg(r,-\ell) \right)$. To see the surjectivity of the map, let us fix $\xi\in \left(\fg(r,\ell)\,+\,\fg(r,-\ell) \right)$ for some integer $\ell$, $0\leq \ell <r$. The formulas in b) are compatible with the Cartan involution and real structure, and consequently determine a unique $\Phi\in\Hom^{\BR,\theta} (\fs,\fg(r))$ such that the formulas in b) hold and $\Phi\,h=\xi$. The criterion (4.14b) and a short computation show that $\Phi$, thus defined, lies in $\Hom^{\BR,\theta}(\fs,\fg(r))(r-2)$. This completes the verification of (8.52).

The definition of the map $f_t$, coupled with the formula for the differential of the exponential map -- see \cite{He,theorem II.1.7}, for example -- leads to the formula
$$
(f_t)_*\zeta \, = \, \Ad \exp(sth)\left(\zeta \ + \ \left[\frac {1\,-\,e^{-st\,\ad h}}{s\,\,t\,\ad h}\, (t \operatorname{Re}\zeta)\,,\, s(h\,+\,ie\,+\,if)\right]\, \right)\,.
\tag8.53
$$
We apply this to $\zeta = \Phi(h+ie+if)$, with $\Phi \in \Hom^{\BR,\theta}(\fs,\fg(r))(r-2)$ viewed as tangent vector to $\Cal O_\fp$ at $\nu_t$, as in (8.50). To simplify the statement we are about to make, we assume
$$
\Phi\,h \ \in \ \left(\fg(r,\ell)\,+\,\fg(r,-\ell) \right) \cap \pr \qquad (\,0\leq \ell < r\,)\,.
\tag8.54
$$
In any case, $\,\Hom^{\BR,\theta}(\fs,\fg(r))(r-2)\,$ has a basis consisting of linear maps $\Phi$ of this type. According to (8.52b), $\Phi e$ has a nonzero component in $\fg(r,\ell+2)$ -- unless $\Phi = 0$, of course -- but no components in $\fg(r,j)$ with $j>\ell+2$; similarly, $\Phi f$ has no components in $\fg(r,j)$ with $j\geq \ell$. The operator $\,\Ad \exp(sth)\,$ acts on $\fg(r,\ell)$ as multiplication by $e^{st\ell}=s^{-\ell/2}$, whereas the operator $\,(\ad h)^{-1}(1\,-\,e^{-st\ad h})\,$ acts by $\ell^{-1}(1-e^{-st\ell})=\ell^{-1}(1-s^{\ell/2})$ or $st = -\frac 12 \log s$, depending on whether $\ell>0$ or $\ell=0$. Looking at the leading terms, or equivalently the terms involving the lowest power of $s$, we find
$$
(f_t)_*|_{\nu_t}\,\Phi(h+ie+if) \ = \ 
\cases
i\,s^{-1- \frac \ell 2}\,\frac {1 \, +\, \ell(r-\ell)}{\ell(r-\ell)}\,(\Phi e)_{\ell+2}\ +\ \dots  &\text{if}\ \ell > 0
\\
i\,|\log s|\,s^{-1}\,\frac 1{2r}\,\Phi e\ \ + \, \dots &\text{if}\ \ell = 0 \ ;
\endcases
\tag8.55
$$
here $\,\dots\,$ refers to lower order terms, and we are using (8.52b) to express $[(\Phi h)_\ell , e]$ as a multiple of $(\Phi e)_{\ell+2}\,$.

Let us re-state the top line of the identity (8.55) in terms $\Phi e$, rather than $(\Phi e)_{\ell + 2}$. By (8.54) and (8.52b), if $\ell > 0$,
$$
\aligned
&\frac {2r}{r+\ell}\,(\Phi e)_{\ell +2} \ = \ \frac {-2r}{(r-\ell)(r+\ell)}\,[\,e\,,\,(\Phi h)_\ell\,]
\\
&\qquad =\ \frac {-1}{r - \ell}\,[\,e\,,\,(\Phi h)_\ell\,]\, - \, \frac 1{r + \ell}\,[\,e\,,\,(\Phi h)_{-\ell}\,]\, - \, \frac 1{r + \ell}\,[\,e\,,\,(\Phi h)_{\ell}\,-\,(\Phi h)_{-\ell}\,]
\\
&\qquad =\ (\Phi e)_{\ell + 2} \ + \ (\Phi e)_{-\ell + 2} \ - \ \frac 1{r+\ell}\,[\,e\,,\,(\Phi h)_{\ell}\,+\,\theta((\Phi h)_\ell)\,]
\\
&\qquad =\ \Phi e \ - \ \frac 1{r+\ell}\,[\,e\,,\,(\Phi h)_{\ell}\,+\,\theta((\Phi h)_{-\ell})\,]\,.
\endaligned
$$
In the next to last line, we have used the fact that $\theta$ acts as $-1$ on $\Phi h$ and maps $\fg(r,\ell)$ to $\fg(r,-\ell)$. Thus
$$
\aligned
&(\Phi e)_{\ell +2} \ = \ \frac {r + \ell}{2r}\,\Phi e \ - \ [\,e\,,\,\eta_\Phi\,]\,,
\\
&\qquad\qquad \text{with}\ \ \eta_\Phi \ =_{\text{def}}\ \frac 1{2r}\,\left(\,(\Phi h)_\ell\,+\,\theta((\Phi h)_\ell)\,\right)\ \in \ \kr\,.
\endaligned
\tag8.56
$$
Combining (8.55-56), we get 
$$
\alignedat2
&{\lim}_{t\to\infty}\,\left( s^{1+\frac\ell 2}\, (f_t)_*|_{\nu_t}\,\Phi(h+ie+if) \right)\ =
\\
&\qquad\qquad \frac {r + \ell + \ell (r^2 - \ell^2)}{2r \ell (r-\ell)}\, \Phi(ie)\ + \ \frac{1 + \ell (r-\ell)}{\ell(r-\ell)}\,[\,\eta_\Phi\,,\,ie\,] \ \ \ &\text{if}\ \ \ell > 0\,,
\\
&{\lim}_{t\to\infty}\,\left( \frac s {|\log s|}\, (f_t)_*|_{\nu_t}\,\Phi(h+ie+if) \right)\ = \ \Phi(ie) \ \ \ &\text{if}\ \ \ell = 0\,.
\endalignedat
\tag8.57
$$
In analogy to (8.50), we can describe the tangent space to $\Cal O_\BR$ at the point $ie$ as
$$
\gathered
T_{ie}\Cal O_\BR \ = \ T_{ie} (\KR\cdot (ie)) \ \tsize\bigoplus \ \fd_\BR(\Phi_0)(ie) \,,\ \ \ \ \text{with}
\\
\fd_\BR(\Phi_0)\ = \ \tsize\bigoplus_{r\geq 2} \Hom^{\BR,\theta}(\fs,\fg(r))(r-2)\,,\ \ \ T_{ie}\Cal O_\BR \ = \ [\,\kr\,,\,ie\,]\,.
\endgathered
\tag8.58
$$
We have established (8.49); equivalently, there exists a basis $\{\eta_j(t)\}$ of $(f_t)_*(T_{\nu_t \Cal O_\fp})$, depending continuously on the parameter $t$, such that the limits $\tilde\eta_j = \lim_{t\to\infty}\eta_j(t)$ exist and constitute a basis of $T_{ie}\Cal O_\BR$. This follows from the analogous statement about the tangent spaces of the $\KR$-orbits -- which is equivalent to (8.51) -- in conjunction with (8.50), (8.52a), (8.57-58), and the non-vanishing of the coefficients of $\Phi(ie)$ in (8.57). We have pointed out already that (8.49) implies the transversality assertion of the lemma.

To pin down the sign of the intersection, it suffices to compare two orientations on $T_{ie}\Cal O_\BR = \lim_{t\to\infty} ((f_t)_* T_{\nu_t} \Cal O_\fp)$ -- on the one hand, the orientation introduced by the symplectic form $\frac 1{2\pi i}\sigma_\Cal O$, on the other, the orientation coming from the complex structure on $\Cal O_\fp \simeq f_t(\Cal O_\fp)$ and the limiting process; the sign of the intersection is the sign which relates the two orientations. We had remarked already that the tangent spaces $T_{\nu_t} \Cal O_\fp$ all coincide when we regard them as subspaces of $\fg^\BR$. In particular, they all coincide with the tangent space at $\nu_0 = h+ie+if$:
$$
T_{\nu_t} \Cal O_\fp\, = \, T_{\nu_0} \Cal O_\fp \, = \, [\fk,\nu_0]\,.
\tag8.59
$$
For reasons of continuity, the real 2-form $\,\operatorname{Im} \sigma_\Cal O$ is non-degenerate on $(f_t)_*(T_{\nu_t} \Cal O_\fp)$ for all large enough values of $t$. We must show that $(f_t)_*$ is orientation preserving with respect to this symplectic structure on $(f_t)_*(T_{\nu_t} \Cal O_\fp)$ and the orientation of $T_{\nu_0} \Cal O_\fp$ as complex vector space -- equivalently, that $\,\operatorname{Im} (f_t^* \sigma_\Cal O)$, for $t\gg 0$, orients the tangent space $T_{\nu_0} \Cal O_\fp = T_{\nu_t}\Cal O$ consistently with the complex orientation. In fact, we shall show
$$
\aligned 
\text{a)}\ \ &\operatorname{Im} \sigma_\Cal O \ \ \text{is non-degenerate on} \ \ (f_t)_*(T_{\nu_t} \Cal O_\fp)\ \ \text{for all}\ \ t>0\,;
\\
\text{b)}\ \ &S \ =_{def} \ {\lim}_{t\to 0^+} \left( t^{-1} \, f_t^* (\operatorname{Im} \sigma_\Cal O)|_{T_{\nu_t}\Cal O} \right) \ \ \text{exists, is non-degenerate,}
\\
&\text{and orients}\ \ T_{\nu_0}\Cal O = T_{\nu_t}\Cal O\ \ \text{consistently with the complex structure}.
\endaligned
\tag8.60
$$
That suffices: the 2-forms $f_t^* (\operatorname{Im} \sigma_\Cal O)$, for $t>0$, are then all nondegenerate on $T_{\nu_0}\Cal O$ and therefore induce the same orientation. Because of b), this orientation agrees with the orientation determined by the complex structure.

We break down the verification of the statement (8.60a) into the following two separate assertions:
$$
\aligned 
\text{a)}\ \ &\text{the submanifolds} \ \ f_t(\Cal O_\fp) \ \ \text{of the complex orbit}\ \ \Cal O \ \ \text{are}
\\
&\text{Lagrangian with respect to the symplectic form}\ \ \operatorname{Re} \sigma_\Cal O\,;
\\
\text{b)}\ \ &(f_t)_*(T_{\nu_0}\Cal O) \cap i\,(f_t)_*(T_{\nu_0}\Cal O)\ = 0 \ \ \text{for all}\ \ t>0\,.
\endaligned
\tag8.61
$$
Let us assume this for the moment. If $\zeta \in (f_t)_*(T_{\nu_t} \Cal O_\fp)$ lies in the radical of the restriction of $\,\operatorname{Im} \sigma_\Cal O$ to $(f_t)_*(T_{\nu_t} \Cal O_\fp)$, (8.61a) allows us to argue
$$
\gathered 
\operatorname{Im} \sigma_\Cal O (\,\zeta\,,\,(f_t)_*(T_{\nu_t} \Cal O_\fp)\,)\ = \ 0\ \ \ \Longrightarrow \ \ \ \sigma_\Cal O (\,\zeta\,,\,(f_t)_*(T_{\nu_t} \Cal O_\fp)\,)\ = \ 0\ \ \ \Longrightarrow
\\
\sigma_\Cal O (\,\zeta\,,\,i\,(f_t)_*(T_{\nu_t} \Cal O_\fp)\,)\, = \, 0\ \ \Longrightarrow\ \ \sigma_\Cal O (\,\zeta\,,\, (f_t)_*(T_{\nu_t} \Cal O_\fp) \tsize\oplus  i\,(f_t)_*(T_{\nu_t} \Cal O_\fp)\,)\, = \, 0\,;
\endgathered
$$
at the second step we are using the complex linearity of $\sigma_\Cal O$. But (8.61b) and a dimension count imply that $(f_t)_*(T_{\nu_t} \Cal O_\fp)$ and $i\,(f_t)_*(T_{\nu_t} \Cal O_\fp)$ span the tangent space of $\Cal O$ at $f_t(\nu_t)$, so $\zeta$ lies in the radical of the holomorphic symplectic form $\sigma_\Cal O$, forcing $\zeta = 0$. Thus (8.61a,b) do imply (8.60a). At this point, only (8.60b) and (8.61a,b) remain to be proved.

Recall the notation (8.47) and the formula (8.53) for the differential of $f_t$. We apply this formula to a tangent vector $\zeta\in T_{\nu_t}\Cal O_\fp$. Because of (8.59), we can write $\zeta = [\kappa,\nu_0]$ for some $\kappa\in\fk$\,, so that
$$
\aligned
&(f_t)_*\zeta\ = \ (f_t)_*[\kappa,\nu_0]\ =\vspace{-2\jot}
\\
&\qquad s^{-1}\,\left[ \Ad\exp(sth)\left( \kappa \,+\, \frac {1 - e^{-st\ad h}}{\ad h}\operatorname{Re}[\kappa,\nu_0]\right)\,,\, \Ad\exp(sth)\,\nu_t\, \right] \ =\vspace{0\jot}
\\
&\qquad -\,s^{-1}\,\left[\,f_t(\nu_t)\,,\, \Ad\exp(sth)\left( \kappa \,+\, \frac {1 - e^{-st\ad h}}{\ad h}\operatorname{Re}\zeta\right) \, \right]\ .
\endaligned
\tag8.62
$$
The holomorphic symplectic form $\sigma_\Cal O$ is the canonical symplectic form of the complex coadjoint orbit that corresponds to $\Cal O$ when we identify $\fg \simeq \fg^*$ via the Killing form. Thus, for  $\zeta_j = [\kappa_j,\nu_0] \in T_{\nu_t}\Cal O_\fp$, $j=1,2$,
$$
\aligned
&(f_t^*\sigma_\Cal O)(\zeta_1,\zeta_2)\ = \ (\sigma_\Cal O)|_{f_t(\nu_t)}(\,(f_t)_*[\kappa_1,\nu_0]\,,\,(f_t)_*[\kappa_2,\nu_0]\,)\ = \vspace{-2\jot}
\\
&B\left(\,f_t(\nu_t)\,,\,\left[\, (\ad f_t(\nu_t))^{-1}((f_t)_*[\kappa_1,\nu_0])\,,\, (\ad f_t(\nu_t))^{-1}((f_t)_*[\kappa_2,\nu_0])\,\right]\,\right)\ =\vspace{0\jot}
\\
&s^{-2}\, B\left(\,\nu_t\,,\,\left[ \, \kappa_1 \,+\, \frac {1 - e^{-st\ad h}}{\ad h}\operatorname{Re}\zeta_1\,,\, \kappa_2 \,+\, \frac {1 - e^{-st\ad h}}{\ad h}\operatorname{Re}\zeta_2\,\right]\,\right)\ =\vspace{0\jot}
\\
&s^{-2}\, B\left(\,\nu_t\,,\,\left[ \, \frac {1 - e^{-st\ad h}}{\ad h}\operatorname{Re}\zeta_1\,,\, \frac {1 - e^{-st\ad h}}{\ad h}\operatorname{Re}\zeta_2\,\right]\,\right)\ +\vspace{0\jot}
\\
&s^{-2}\, B\left(\,\nu_t\,,\,\left[ \, \kappa_1 \,,\, \frac {1 - e^{-st\ad h}}{\ad h}\operatorname{Re}\zeta_2\,\right] \ + \ \left[ \, \frac {1 - e^{-st\ad h}}{\ad h}\operatorname{Re}\zeta_1\,,\, \kappa_2\,\right]\,\right)\,.
\endaligned
$$
Here, in the second line, $(\ad f_t(\nu_t))^{-1}((f_t)_*[\kappa_j,\nu_0])$ is symbolic notation for any element of $\fg$ whose image under $\,\ad f_t(\nu_t)$ is $(f_t)_*[\kappa_j,\nu_0]$; in passing from the second line to the third, we are using (8.62), the identity $f_t(\nu_t)= \Ad\exp(st\ad h)(\nu_t)$, and the $\Ad$-invariance of the Killing form; the last step is justified by the perpendicularity of $\nu_t\in\fp$ and $[\kappa_1,\kappa_2] \in \fk$. Next, we use the infinitesimal invariance of $\,B\,$, the relation $\nu_t = s\,\nu_0$, and the relation between $\zeta_j$ and $\kappa_j$, to conclude
$$
\aligned
&(f_t^*\sigma_\Cal O)(\zeta_1,\zeta_2)\ = \ (\sigma_\Cal O)|_{f_t(\nu_t)}(\,(f_t)_*[\kappa_1,\nu_0]\,,\,(f_t)_*[\kappa_2,\nu_0]\,)\ = \vspace{-2\jot}
\\
&=\ s^{-1}\, B\left(\,\nu_0\,,\,\left[ \, \frac {1 - e^{-st\ad h}}{\ad h}\operatorname{Re}\zeta_1\,,\, \frac {1 - e^{-st\ad h}}{\ad h}\operatorname{Re}\zeta_2\,\right]\,\right)\vspace{0\jot}
\\
&-\,s^{-1}\, B\left(\,\zeta_1 \,,\, \frac {1 - e^{-st\ad h}}{\ad h}\operatorname{Re}\zeta_2\, \right)\  + \ s^{-1}\, B\left( \, \frac {1 - e^{-st\ad h}}{\ad h}\operatorname{Re}\zeta_1\,,\, \zeta_2\,\right)\,.
\endaligned
\tag8.63
$$
We shall use this formula to verify (8.60b) and (8.61a).

For $t$ near 0, $s(t)=1-2t+\dots\,$ and $\,(\ad h)^{-1}(1 - e^{-st\ad h}) = st\!\cdot \!1 + \dots\,$, hence $\,s^{-1}(\ad h)^{-1}(1 - e^{-st\ad h}) = t\!\cdot \!1 + \dots\,$, and
$$
f_t^*(\operatorname{Im}\sigma_\Cal O)(\zeta_1,\zeta_2)\ = \ -\,t\,B(\operatorname{Im}\zeta_1\,,\,\operatorname{Re}\zeta_2)\,+ \,t\,B(\operatorname{Re}\zeta_1\,,\,\operatorname{Im}\zeta_2)\,+\,\dots\ .
\tag8.64
$$
We conclude that $S = \lim_{t\to 0^+} t^{-1} f_t^*(\operatorname{Im}\sigma_\Cal O)$ exists as $\BR$-bilinear, alternating form on $T_{\nu_0}\Cal O_\fp = [\fk,\nu_0]$ and is given by the formula
$$
S(\zeta_1,\zeta_2) \ = \ -\, B(\operatorname{Im}\zeta_1\,,\,\operatorname{Re}\zeta_2)\,+ \,B(\operatorname{Re}\zeta_1\,,\,\operatorname{Im}\zeta_2)\,.
\tag8.65
$$
Let $\{\zeta_j\}$ be a $\BC$-basis of $[\fk,\nu_0]$, orthonormal with respect to the inner product (2.2). Since $\theta$ acts as multiplication by $-1$ on $[\fk,\nu_0] \subset \fp$, 
$$
\gathered
S(\zeta_j,i\,\zeta_k)\, = \,  B(\operatorname{Im}\zeta_j\,,\,\operatorname{Im}\zeta_k)+ B(\operatorname{Re}\zeta_j\,,\,\operatorname{Re}\zeta_k)\, = \, \operatorname{Re}(\zeta_j,\zeta_k)\,=\,\delta_{j,k}\,,
\\
S(\zeta_j,\zeta_k)\, = \,  -\, B(\operatorname{Im}\zeta_j\,,\,\operatorname{Re}\zeta_k)+ B(\operatorname{Re}\zeta_j\,,\,\operatorname{Im}\zeta_k)\, = \, -\,\operatorname{Im}(\zeta_j,\zeta_k)\,=\,0\,.
\endgathered
\tag8.66
$$
In particular, the nondegenerate alternating bilinear form $S$ orients $[\fk,\nu_0]$, viewed as real vector space, in the same way as the complex structure. This establishes (8.60b).

The formula (8.63) and its derivation remain valid if we replace $\nu_t = s\,\nu_0$ by an arbitrary point $\nu \in \Cal O_\fp$ and $s\,h = \operatorname{Re}\nu_t$ by $\,\operatorname{Re}\nu$. We take real parts on both sides, to conclude
$$
\aligned
&(\operatorname{Re} \sigma_\Cal O)|_{f_t(\nu)}(\,(f_t)_*\zeta_1\,,\,(f_t)_*\zeta_2\,)\ = \vspace{-2\jot}
\\
&= \ B\!\left(\,\operatorname{Re}\nu\,,\,\left[ \, \frac {1 - e^{-t\ad \, \operatorname{Re}\nu}}{\ad \, \operatorname{Re}\nu}\operatorname{Re}\zeta_1\,,\, \frac {1 - e^{-t\ad \, \operatorname{Re}\nu}}{\ad \, \operatorname{Re}\nu}\operatorname{Re}\zeta_2\,\right]\,\right) \vspace{0\jot}
\\
&-\,B\!\left(\,\operatorname{Re}\zeta_1 \,,\, \frac {1 - e^{-t\ad \, \operatorname{Re}\nu}}{\ad \, \operatorname{Re}\nu}\operatorname{Re}\zeta_2\, \right)\  + \ B\!\left( \, \frac {1 - e^{-t\ad \, \operatorname{Re}\nu}}{\ad \, \operatorname{Re}\nu}\operatorname{Re}\zeta_1\,,\, \operatorname{Re}\zeta_2\,\right)\,,
\endaligned
\tag8.67
$$
for all $\zeta_1,\zeta_2 \in T_\nu \Cal O_\fp$. On the other hand, because of the invariance of $B$,
$$
\aligned
&B\!\left(\,\operatorname{Re}\nu\,,\,\left[ \, \frac {1 - e^{-t\ad \, \operatorname{Re}\nu}}{\ad \, \operatorname{Re}\nu}\operatorname{Re}\zeta_1\,,\, \frac {1 - e^{-t\ad \, \operatorname{Re}\nu}}{\ad \, \operatorname{Re}\nu}\operatorname{Re}\zeta_2\,\right]\,\right)\ =\vspace{0\jot}
\\
&B\!\left(\ad \, \operatorname{Re}\nu\circ \frac {1 - e^{-t\ad \, \operatorname{Re}\nu}}{\ad \, \operatorname{Re}\nu}\operatorname{Re}\zeta_1\,,\, \frac {1 - e^{-t\ad \, \operatorname{Re}\nu}}{\ad \, \operatorname{Re}\nu}\operatorname{Re}\zeta_2\right)\ =\vspace{0\jot}
\\
&B\!\left( ( 1 - e^{-t\ad \, \operatorname{Re}\nu})\operatorname{Re}\zeta_1\,,\, \frac {1 - e^{-t\ad \, \operatorname{Re}\nu}}{\ad \, \operatorname{Re}\nu}\operatorname{Re}\zeta_2\right)\ =\vspace{0\jot}
\\
&B\!\left(\operatorname{Re}\zeta_1\,,\, \frac {1 - e^{-t\ad \, \operatorname{Re}\nu}}{\ad \, \operatorname{Re}\nu}\operatorname{Re}\zeta_2\right)\, - \, B\!\left( \operatorname{Re}\zeta_1\,,\, \frac {e^{t\ad \, \operatorname{Re}\nu}-1}{\ad \, \operatorname{Re}\nu}\operatorname{Re}\zeta_2\right)
\endaligned
\tag8.68
$$
The operator $\ad \, \operatorname{Re}\nu$ is skew with respect to $B$, so $\,(\ad \, \operatorname{Re}\nu)^{-1}(1 - e^{-t\ad \, \operatorname{Re}\nu})$ is the adjoint of $\,(\ad \, \operatorname{Re}\nu)^{-1}(e^{t\ad \, \operatorname{Re}\nu}-1)$, and
$$
B\!\left( \operatorname{Re}\zeta_1\,,\, \frac {e^{t\ad \, \operatorname{Re}\nu}-1}{\ad \, \operatorname{Re}\nu}\operatorname{Re}\zeta_2\right)\ = \ B\!\left( \frac {1 - e^{-t\ad \, \operatorname{Re}\nu}}{\ad \, \operatorname{Re}\nu}\operatorname{Re}\zeta_1\,,\, \operatorname{Re}\zeta_2\right)\,.
\tag8.69
$$
Combining (8.67-69), we find that $\,\operatorname{Re}\sigma_\Cal O$ vanishes identically on $f_t(\Cal O_\fp)$. Since $\Cal O_\fp$ has half the dimension of $\Cal O$, this implies (8.61a).

Only (8.61b) remains to be established. Let us assume, then, that $t>0$. We consider two tangent vectors
$$
\zeta_j \in T_{\nu_t}\Cal O_\fp = [\,\fk\,, h+ei+if]\subset \fp\,,\ \ \text{such that}\ \ (f_t)_*\zeta_1\ = \ i\,(f_t)_*\zeta_2\,.
\tag8.70
$$
We express the $\zeta_j$ in terms of their real and imaginary parts,
$$
\zeta_j\ = \ \xi_j\,+\,i\,\eta_j\,,\qquad \text{with}\ \ \xi_j\,,\, \eta_j \in \pr\,.
\tag8.71
$$
Because $\,(\ad h)^{-1}(1-e^{-st \ad h})[\xi_j,h] = (e^{-st \ad h}-1)\xi_j$, the formula (8.53) can be re-written as follows:
$$
(f_t)_*\zeta_j\ = \ \xi_j\ + \ i\,\Ad\exp(sth)\left(\eta_j\,+\, \left[ \frac{1- e^{-st\ad h}}{\ad h}\,\xi_j\,,\,e+f\right] \right)\,.
\tag8.72
$$
Our  assumption (8.70) on the $\zeta_j$ is therefore equivalent to
$$
\aligned
e^{-s\,t\,\ad\,h}\,\xi_1\ &= \ -\,\eta_2\,+\, \left[e+f \,,\, \frac{1- e^{-st\ad h}}{\ad h}\,\xi_2\right]\,,
\\
e^{-s\,t\,\ad\,h}\,\xi_2\ &= \ \eta_1\,-\, \left[e+f \,,\,\frac{1- e^{-st\ad h}}{\ad h}\,\xi_1\right]\,.
\endaligned
\tag8.73
$$
We need to separate the components in $\kr$ and $\pr$. For this purpose, we define
$$
\aligned
S \ &= \ {\operatorname{sinh}(s\,t\,\ad\,h)\over\ad \, h} \ = \ s\,t\cdot 1\ + \ {1 \over 6}(s\,t\,\ad\,h)^2\ + \ \dots\ ,
\\
T \ &= \ {1\,-\,\operatorname{cosh}(s\,t\,\ad\,h)\over\ad \, h} \ = \ -\,\frac 12 \,s\,t\,\,\ad\,h\ - \ {1 \over 24}(s\,t\,\ad\,h)^3\ - \ \dots\ .
\endaligned
\tag8.74
$$
Even powers of $\,\ad\, h\,$ or $\,\ad(e+f)\,$ commute with the Cartan involution, whereas odd powers anti-commute; also,
$$
e^{-s\,t\,\ad\,h}\ = \ 1\ -\ (\,S \, +\, T\,)\circ \ad\,h\ = \ 1\ - \ \ad \,h \circ (\,S \, +\, T\,)\,.
\tag8.75
$$
Equating $\pr$-components in (8.73), we now find
$$
\aligned
\eta_1 \ &= \ (1\,-\,\ad\,h \circ T) \xi_2\ + \ \ad(e+f)\circ T\ \xi_1\,,
\\
\eta_2 \ &= \ -(1\,-\,\ad\,h \circ T) \xi_1\ + \ \ad(e+f)\circ T\ \xi_2\,,
\endaligned
\tag8.76a
$$
and the equality of the $\kr$-components translates into
$$
\aligned
\ad\,h \circ S\ \xi_1 \ &= \ -\,\ad(e+f)\circ S\ \xi_2\,,
\\
\ad\,h \circ S\ \xi_2 \ &= \ \ad(e+f)\circ S\ \xi_1\,.
\endaligned
\tag8.76b
$$
The latter two equations can be combined into the single complex equation
$$
[\,h\,+\,i\,e\,+\,i\,f\,, \, S(\xi_1\,-\,i\,\xi_2)\,] \ = \ 0\,.
\tag8.77
$$
We shall use these equations to show that the $\zeta_j$ must vanish..

Both $\zeta_j$ lie in $[\fk, h+ie+if]$, hence in the image of $\,\ad(h+ie+if):\fg\to\fg$\,, and $S(\xi_1-\xi_2)$ lies in the kernel of $\ad(h+ie+if)$ by (8.77). The image and the kernel are each other's annihilator, relative to the Killing form. Thus
$$
B(\,S(\xi_1-i\xi_2)\,,\,\zeta_1 + i \zeta_2\,)\ = \ 0\,.
$$
Taking real parts, we find
$$
\aligned
0\ &= \ B(S\xi_1,\xi_1-\eta_2)\, +\, B(S\xi_2,\xi_2 + \eta_1) 
\\
&= \ B(S\xi_1,\xi_1)\, +\, B(S\xi_1,(1-\ad h \circ T)\xi_1)\, -\, B(S\xi_1,\ad(e+f) \circ T\xi_2)
\\
&\ \ \ \ +\, B(S\xi_2,\xi_2)\, +\, B(S\xi_2,(1-\ad h \circ T)\xi_2)\, +\, B(S\xi_2,\ad(e+f) \circ T\xi_1)\,;
\endaligned
\tag8.78
$$
at the second step, we have used (8.76a) to express the $\eta_j$ in terms of the $\xi_j$. The infinitesimal invariance of the Killing form and (8.76b) give
$$
\aligned
&B(S\xi_1,\ad(e+f) \circ T\xi_2)\ = \ -\,B(\ad(e+f) \circ S\xi_1,T\xi_2)\ =  
\\
&\qquad -\,B(\ad h \circ S\xi_2,T\xi_2)\ = \ B(S\xi_2,\ad h \circ T\xi_2)\,,
\endaligned
\tag8.79a
$$
and similarly
$$
\aligned
&B(S\xi_2,\ad(e+f) \circ T\xi_1)\ = \ -\,B(\ad(e+f) \circ S\xi_2,T\xi_1)\ =  
\\
&\qquad B(\ad h \circ S\xi_1,T\xi_1)\ = \ -\,B(S\xi_1,\ad h \circ T\xi_1)\,.
\endaligned
\tag8.79b
$$
The operators
$$
1\,-\,\ad\, h \circ T \ = \ \cosh(s\,t\,\ad\,h)\,,\ \ \ \ad\, h \circ T \ = \ 1\,-\,\cosh(s\,t\,\ad\,h)
\tag8.80
$$
are series in $(\ad h)^2$, hence symmetric with respect to the Killing form. Thus, combining (8.78-90), we find
$$
\aligned
0\ &= B(S\xi_1,\xi_1)\, +\, B((1-\ad h \circ T)\circ S\xi_1,\xi_1)\, -\, B(\ad h \circ T \circ S\xi_2,\xi_2)
\\
&+\, B(S\xi_2,\xi_2)\, +\, B((1-\ad h \circ T)\circ S\xi_2,\xi_2)\, -\, B(\ad h \circ T \circ S\xi_1,\xi_1)\,.
\endaligned
\tag8.81
$$
The inner product (2.2) agrees with the Killing form on $\pr$. Relative to this inner product, $\ad h$ is a symmetric operator, whose eigenspace decomposition diagonalizes $S$ and $\,\ad h \!\circ \! T$. For $t>0\,$ -- \,which also makes $s$ strictly positive -- the eigenvalues of $S$ and $1-\ad h \!\circ \! T$ are strictly positive, and those of $\,\ad h \!\circ \! T$ non-positive. Thus all terms in (8.81) vanish individually, and $\xi_1 = \xi_2 = 0$. The $\eta_j$, which can be expressed in terms of the $\xi_j$, must vanish also. We have shown that (8.70) forces $\zeta_1 = \zeta_2 = 0$. This completes the verification of (8.61b), and with it, the proof of lemma 8.10.
\enddemo

\vskip 3\jot

\enddocument